%% file: main-revision.tex
\documentclass[sn-mathphys-num]{sn-jnl}
\usepackage{graphicx} %
\usepackage{caption,subcaption}
\usepackage[dvipsnames]{xcolor}
\usepackage{listings}
\usepackage{amsmath,amsfonts,amssymb,amsthm}
\usepackage{nicefrac}
\usepackage{todonotes}
\usepackage{multirow}%
\usepackage{mathrsfs}%
\usepackage[title]{appendix}%
\usepackage{textcomp}%
\usepackage{manyfoot}%
\usepackage{booktabs}%
\usepackage{listings}%
\usepackage{caption}
\usepackage{siunitx}

\usepackage{algorithm2e}
\theoremstyle{plain}

\newtheorem{theorem}{Theorem}
\theoremstyle{remark}
\newtheorem{remark}{Remark}
\newtheorem{example}{Example}
\usepackage{booktabs}

\usepackage{tikz}
\usepackage{pgfplots}
\pgfplotsset{compat=1.18}

\usepackage{hyperref}
\usepackage{microtype}

\lstdefinelanguage{julia}
{
  keywordsprefix=\@,
  morekeywords={
    exit,whos,edit,load,is,isa,isequal,typeof,tuple,ntuple,uid,hash,finalizer,convert,promote,
    subtype,typemin,typemax,realmin,realmax,sizeof,eps,promote_type,method_exists,applicable,
    invoke,dlopen,dlsym,system,error,throw,assert,new,Inf,Nan,pi,im,begin,while,for,in,return,
    break,continue,macro,quote,let,if,elseif,else,try,catch,end,bitstype,ccall,do,using,module,
    import,export,importall,baremodule,immutable,local,global,const,Bool,Int,Int8,Int16,Int32,
    Int64,Uint,Uint8,Uint16,Uint32,Uint64,Float32,Float64,Complex64,Complex128,Any,Nothing,None,
    function,type,typealias,abstract
  },
  sensitive=true,
  morecomment=[l]{\#},
  morestring=[b]',
  morestring=[b]" 
}

\begin{document}

\title{Stage-Parallel Implicit Runge--Kutta methods via low-rank matrix equation corrections}

\author[1]{\fnm{Fabio} \sur{Durastante}}\email{fabio.durastante@unipi.it}
\equalcont{These authors contributed equally to this work.}

\author*[2]{\fnm{Mariarosa} \sur{Mazza}}\email{mariarosa.mazza@uniroma2.it}
\equalcont{These authors contributed equally to this work.}

\affil[1]{\orgdiv{Department of Mathematics}, \orgname{University of Pisa}, \orgaddress{\street{Largo Bruno Pontecorvo, 5}, \city{Pisa}, \postcode{56127}, \state{PI}, \country{Italy}}}

\affil*[2]{\orgdiv{Department of Mathematics}, \orgname{University of Roma ``Tor Vergata''}, \orgaddress{\street{Via Della Ricerca Scientifica, 1}, \city{Rome}, \postcode{00133}, \state{RM}, \country{Italy}}}

\abstract{
Implicit Runge--Kutta (IRK) methods are highly effective for solving stiff ordinary differential equations (ODEs) but can be computationally expensive for large-scale problems due to the need of solving coupled algebraic equations at each step. This study improves IRK efficiency by leveraging parallelism to decouple stage computations and reduce communication overhead, specifically we stably decouple a perturbed version of the stage system of equations and recover the exact solution by solving a Sylvester matrix equation with an explicitly known low-rank right-hand side. Two IRK families---symmetric methods and collocation methods---are analyzed, with extensions to nonlinear problems using a simplified Newton method. Implementation details, shared memory parallel code, and numerical examples, particularly for ODEs from spatially discretized PDEs, demonstrate the efficiency of the proposed IRK technique.}

\keywords{Runge--Kutta; Stage-Parallel Method; Low-rank Sylvester Matrix Equation.}

\maketitle

\section{Introduction}

Implicit Runge--Kutta (IRK) methods are a widely used technique for the numerical integration of ordinary differential equations (ODEs):  
\begin{equation}\label{eq:theproblem}
\begin{cases}
    M \mathbf{y}'(t)=f(\mathbf{y}(t),t), & t \in [0,T],  \\
    \mathbf{y}(0) = \mathbf{y}_0,
\end{cases} \quad \begin{array}{l}
\mathbf{y} \,:[0,T] \rightarrow \mathbb{R}^N,\\
f \,:\,\mathbb{R}^N \times [0,T] \rightarrow \mathbb{R}^N, \\
\mathbf{y}_0 \in \mathbb{R}^N, \, M \in \mathbb{R}^{N \times N},\; N \in \mathbb{N}, \, T > 0.
\end{array}
\end{equation}
Unlike explicit methods, which compute the solution at the next time-step using only information from previous stages, an IRK method with $s$ stages involves solving a system of algebraic equations at each step:
\begin{equation}\label{eq:steps_Runge_Kutta}
M \mathbf{y}^{(n+1)} = M \mathbf{y}^{(n)} + h \sum_{i=1}^{s} b_i \mathbf{k}_{i}^{(n)}, \quad n = 0, \ldots, n_t, \quad h = \frac{T}{n_t}, \quad n_t \in \mathbb{N},
\end{equation}
where the stages $\mathbf{k}_{i}^{(n)}$ for the approximation $\mathbf{y}^{(n+1)}\approx \mathbf{y}(t_{n+1})$ are given by
\begin{equation}\label{eq:stages_Runge_Kutta}
\mathbf{k}_{i}^{(n)} = f\left(\mathbf{y}^{(n)} + h \sum_{j=1}^{s} a_{i,j} \mathbf{k}_{j}^{(n)}, t_n + c_i \tau\right), \qquad i=1,\ldots,s,
\end{equation}
with $t_n = n h$. The Runge--Kutta method is fully defined by the coefficients $a_{i,j}$, weights $b_i$, and nodes $c_i$, for $i,j = 1, \ldots, s$.  

IRK methods are particularly suitable for stiff ODEs, as explicit methods can become unstable or require impractically small time-steps. By considering the interdependence of intermediate solution values within each time-step, IRK methods offer enhanced stability and accuracy. However, these benefits come at the cost of solving a computationally expensive system of algebraic equations at each step. This challenge is especially pronounced when IRK methods are applied to semi-discretized partial differential equations (PDEs), where the resulting systems can be of large scale.  

To address this computational burden, IRK methods can benefit significantly from parallelism. In particular, parallelization across the stages enables simultaneous computation of multiple intermediate solutions, reducing overall computational time. This approach has been explored in the literature for various classes of RK methods, see, e.g.,~\cite{Bellen,Bellen2,Houwen1,Houwen2,Houwen3,Iserles,MR1790038,MR4735250}. For instance,~\cite{Houwen1,Houwen2,Houwen3} propose a diagonal iteration method to solve the stage equations~\eqref{eq:stages_Runge_Kutta}, which facilitates parallelization through the solution of block-diagonal systems. Similarly,~\cite{Iserles} investigates the sparsity structure of the stage matrix $A = (a_{i,j})_{i,j}$ in~\eqref{eq:stages_Runge_Kutta} using concepts from sparse linear system solvers, decoupling stage dependencies as much as possible. For parallel explicit methods, including Parallel Explicit RK and RK-Nyström methods, we refer to~\cite{Sommeijer1,Cong1,Cong2}. Parallelism can also be leveraged by employing preconditioned distributed iterative solvers to independently solve the equations arising at each stage~\cite{SouthworthI,SouthworthII,MR1790038,SantoloSVD,MR4735250}. 

In this work, we aim to solve system~\eqref{eq:stages_Runge_Kutta} via the following two-step parallel approach:  
\begin{enumerate}
    \item  \textbf{Perturbation and Decoupling:} In the first step, we perturb the coefficients $a_{i,j}$ in~\eqref{eq:stages_Runge_Kutta}, decoupling the solution. This perturbation enables parallel computation of a modified set of stages. Of course, the resulting solution is affected by an error due to the perturbation. 
    \item \textbf{Correction via Sylvester Equation:} To recover the exact stages, we solve a Sylvester matrix equation with a known low-rank term~\cite{SimonciniReview}. This sequential step uses scalable and efficient Krylov projection methods.  
\end{enumerate}
Once these two steps are completed, we compute the approximation of the ODE solution at the new time-step using~\eqref{eq:steps_Runge_Kutta}.  The concept of splitting the procedure into two steps is inspired by the approach in~\cite{MR4646959}, later expanded in~\cite{MR4713232}, which were developed for constructing parallel-in-time methods for multistep linear methods. In this work, however, instead of using the perturbation to decouple the problem via the Fast Fourier Transform, we leverage the hidden structural properties of the stage matrix $A$. The resulting approach is particularly effective because the number of stages is independent of the problem size, {hence we achieve the decoupling-diagonalization by computing the spectral decomposition of an $s \times s$ matrix in $s^3$ operations. This cost is smaller than either $n_t \log_2(n_t)$, incurred in using FFT-based strategy, or the $n_t^3$ cost of a direct approach.} %
{For all the three choices the time-stepping cost needs to be added}. The construction supports two types of parallelism: distributed memory parallelism, implemented using MPI subcommunicators, and shared memory parallelism. In this discussion, we will focus on the shared memory model of computation.

The remainder of this article is structured as follows: In Section~\ref{sec:parallelism-across-the-stages}, we introduce the general framework for a generic IRK method applied to linear problems. We then analyze and construct two prominent families of IRK methods: symmetric methods (Section~\ref{sec:symmetric-RK}), collocation methods (Section~\ref{sec:collocation_methods}). In Section~\ref{sec:nonlinear}, we extend our analysis to nonlinear problems using a simplified Newton method for~\eqref{eq:theproblem}. Section~\ref{sec:code_and_implementation} details the implementation and accompanying code and provides numerical examples, including ODEs arising from spatially discretized PDEs. Section~\ref{sec:conclusions} concludes the work with suggestions for future research.

\paragraph{Notation} We will denote vectors in bold face, matrices with capital letters (either latin, greek or calligraphic). The $\|\cdot\|$ denotes the two-norm, while $\|\cdot\|_F$ denotes the Frobenius norm. The vectors of the canonical basis will be denoted by $\mathbf{e}_j$, while $\mathbf{1}_s$ denotes the {column} vector of all ones of size $s$. The identity matrix of dimension $m$ is denoted by $I_m$. We call $\operatorname{vec}\left(\cdot\right)$ the operator stacking the columns of its matrix argument one below the other. The symbol ``$\otimes$'' denotes the Kronecker product.

\section{Parallelism across the stages for linear problems}\label{sec:parallelism-across-the-stages}

To succinctly represent the generic IRK method \eqref{eq:steps_Runge_Kutta}-\eqref{eq:stages_Runge_Kutta} and the parallelism across the stage idea we employ the Butcher tableau:
\begin{equation}\label{eq:butcher_tableau}
\def\arraystretch{1.2}
\begin{array}[b]{c|c}
\mathbf{c} & A\\
\hline
 & \mathbf{b}^\top
\end{array} \;\raisebox{0.8em}{=}\; 
\begin{array}[b]{c|ccc}
c_1 & a_{1,1} & \ldots & a_{1,s}\\
\vdots & \vdots & \ddots & \vdots\\
c_s & a_{s,1} & \ldots & a_{s,s}\\
\hline
 & b_1 & \ldots & b_{s}
\end{array}\raisebox{0.8em}{.}
\end{equation}
To describe the stage-parallel version, we consider a case in which the ODE system is given by the semi-discretization in space of a PDE. We examine a linear PDE of the form
\begin{equation}\label{eq:linear-pde}
\begin{cases}
\frac{\partial u}{\partial t} + \mathcal{L}[u] = {s}(\mathbf{x},t), & (\mathbf{x},t) \in \Omega \times (0,T] , \quad \Omega \subseteq \mathbb{R}^{d},   \\
u(0) = u_0(\mathbf{x}) , \\
\mathcal{B}\left[ u(\mathbf{x},t) \right] = 0, &  (\mathbf{x},t) \in \partial \Omega \times [0,T],
\end{cases},\quad d \in \mathbb{N},
\end{equation}
where $\mathcal{L}$ and $\mathcal{B}$ are, respectively, a differential operator involving only spatial derivatives, and a boundary operator.

After a discretization scheme is applied to~\eqref{eq:linear-pde} we can rewrite it in the form~\eqref{eq:theproblem}~as
\begin{equation}\label{eq:linear_problem}
    \begin{cases}
M \mathbf{y}'(t) = - L \mathbf{y}(t) + \hat{\mathbf{f}}(t), \\
\mathbf{y}(0) = \mathbf{y}_0,
\end{cases} \qquad M,L \in \mathbb{R}^{N \times N},
\end{equation}
where the matrices $M$, $L$ and the vector function $\hat{\mathbf{f}}$ collect the discretization of the operators $\mathcal{L}$, $\mathcal{B}$, and of the source term. Applying the generic IRK method with $s$ stages in~\eqref{eq:butcher_tableau} yields
	$$
	\begin{cases}
		\displaystyle M \mathbf y^{(n+1)} = M \mathbf y^{(n)}+h\sum_{i=1}^sb_i \mathbf k_i^{(n)},&n=0,\dots n_t,\\
		\displaystyle M \mathbf k_i^{(n)}= - L \mathbf y^{(n)} -h\sum_{j=1}^s a_{i,j} L \mathbf k_j^{(n)} + \hat {\mathbf{f}}(t_n+c_ih),&i=1,\dots,s,
	\end{cases}
	$$
	where the coefficients $b_i$ and $c_i$ are the nodes and the weights of the Butcher tableau~\eqref{eq:butcher_tableau}. By considering the matrix formulation and setting \[K^{(n)}=[\mathbf k^{(n)}_1|\dots| \mathbf k^{(n)}_s]\in\mathbb R^{N\times s},\] we can rewrite the second equation as
	\begin{align}\label{eq:rk1}
            M K^{(n)}= - L \mathbf{y}^{(n)} \mathbf{1}_s^\top - h L K^{(n)} A^\top + F^{(n)},
	\end{align}
	where
	\[F^{(n)}=[\hat{\mathbf{f}}(t_n+c_1h)|\dots|\hat{\mathbf{f}}(t_n+c_sh)] \in \mathbb{R}^{N \times s}.\]
	If we rewrite~\eqref{eq:rk1} in a linear system format we find
    \begin{equation}\label{eq:block_system}
        	(I_s \otimes M + h A \otimes L ) \operatorname{vec} \left( K^{(n)} \right) = - ( \mathbf{1}_s \otimes L ) \mathbf{y}^{(n)} + \operatorname{vec}\left( F^{(n)}\right).
    \end{equation}
	If we now assume that $A$ is diagonalizable with $A = X \Lambda X^{-1}$ we restate the system as
	\begin{equation*}
	\left(X \otimes I_N\right) \left(I_s \otimes M + h \Lambda \otimes L\right) \left(X^{-1} \otimes I_N\right) \operatorname{vec} \left( K^{(n)} \right) = - ( \mathbf{1}_s \otimes L ) \mathbf{y}^{(n)} + \operatorname{vec}\left( F^{(n)}\right),
	\end{equation*}
	which can then be solved~\cite[\S4.3]{SimonciniReview} as
	\begin{equation}\label{eq:diag-parallel-solve}
	\begin{cases}
	    \mathbf{r} = - ( X^{-1} \mathbf{1}_s \otimes L ) \mathbf{y}^{(n)} + \left(X^{-1} \otimes I_N\right)\operatorname{vec}\left( F^{(n)}\right),\\
	    \left(I_s \otimes M + h \Lambda \otimes L\right) \mathbf{z} = \mathbf{r},\\
	    \operatorname{vec} \left( K^{(n)} \right) =  \left(X \otimes I_N\right) \mathbf{z}.
	\end{cases}
	\end{equation}
    This requires the solution of decoupled linear systems to advance the method in the time-step, i.e., we first compute $\mathbf{r}$
    and then solve the systems
    \begin{equation}\label{eq:linear_systems}
            \left[ M + h \lambda_j L \right] \mathbf{z}^{(j)} = [\mathbf{r}]_{(j-1)N+1:jN}, \qquad j=1,\ldots,s, 
    \end{equation}
    where $[\,\cdot\,]_{p:q}$ extracts the entries from the $p$th to the $q$th, while finally computing the stages as
    \[
    K^{(n)} = Z X^\top, \qquad Z = \left[ \mathbf{z}^{(1)}|\cdots|\mathbf{z}^{(s)}\right].
    \]
    To write this version we assumed that the Butcher Tableau matrix $A$ was diagonalizable. But when is this the case, particularly when is the eigenvector matrix $X$ well-conditioned? In general, this is not the case, and this has led to the development of several types of preconditioners for the system of stages~\eqref{eq:block_system} that exploit the real Schur decomposition~\cite{SouthworthI,SouthworthII} or the SVD factorization~\cite{SantoloSVD} of the matrix $A$.

    In the next two Sections~\ref{sec:symmetric-RK} and~\ref{sec:collocation_methods}, we will show
    how it is possible to perturb the matrix $A$ via a low-rank update to get a diagonalizable matrix for symmetric and collocation IRK methods, respectively. In both cases the solution of the original problem from the perturbed one will be recovered by solving a Sylvester matrix equation.

\section{Symmetric Runge--Kutta schemes}\label{sec:symmetric-RK} An RK scheme is called \emph{symmetric} if it remains invariant under the reflection $t$ to $1-t$ in the direction of integration~\cite{SchererTurke}, i.e., employing the notation in~\eqref{eq:butcher_tableau}, if it is such that
    \[
    \mathbf{c} + J \mathbf{c} = \mathbf{1}_s, \; \mathbf{b} = J \mathbf{b}, \; J A J + A = \mathbf{1}_s\mathbf{b}^\top,
    \]
    for $J$ the exchange matrix{, i.e., the matrix $J$ such that $J [v_1,\ldots,v_{s-1},v_s]^\top = [v_s,v_{s-1},\ldots,v_1]^\top$}; see \cite[\S 5]{BibleVolIII} for further information. 
    
    \begin{example}\label{example:obtaining-coeffs}
    Here we build an example of symmetric IRK scheme. To this end we take $\mathbf{b}$, and the $\mathbf{c}$ as the weights and nodes of the Gauss-Legendre quadrature formula, while the $a_{i,j}$ need to satisfy the sum and moment conditions~\cite[\S\,II.7 p.208]{BibleVol1} and can be computed as integrals of the Lagrange polynomials 
    \begin{equation}\label{eq:gauss2s-A-matrix}
        a_{i,j} = \int_{0}^{c_i} \ell_j(c)\,\mathrm{d}c, \qquad \ell_j(c) = \prod_{\substack{k = 1\\k\neq j}}^{s} \frac{c - c_k}{c_j - c_k}.
    \end{equation}
    The schemes for the low order values of $s$ are given in Table~\ref{tab:symm_rk}, while for the scheme of a higher order the accompanying code\footnote{See the functions in the \texttt{schemes.jl} source code in the GitHub repository \href{https://github.com/Cirdans-Home/SP\_IRK.jl}{Cirdans-Home/SP\_IRK.jl}.} can be employed.
    \begin{table}[htbp]
        \centering
        \caption{Symmetric Runge--Kutta schemes with $s$ stages and order $2s$.}
        \label{tab:symm_rk}
        \begin{tabular}{cc}
        $s = 1$ & $s = 2$ \\
        $ \begin{array}{c|c}
           \nicefrac{1}{2}  & \nicefrac{1}{2}  \\
            \hline
             & 1
        \end{array} $     &  
        $ \begin{array}{c|cc}
           \nicefrac{3-\sqrt{3}}{6}  & \nicefrac{1}{4} & \nicefrac{(3-2\sqrt{3})}{12}  \\
           \nicefrac{3+\sqrt{3}}{6}  & \nicefrac{(3+2\sqrt{3})}{12} & \nicefrac{1}{4}   \\
            \hline
             & \nicefrac{1}{2} & \nicefrac{1}{2}
        \end{array} $   \\
        \midrule
        \multicolumn{2}{c}{$s=3$} \\
        \multicolumn{2}{c}{
        $ \begin{array}{c|ccc}
           \nicefrac{5-\sqrt{15}}{10}  & \nicefrac{5}{36} & \nicefrac{(10-3\sqrt{15})}{45} & \nicefrac{(25-6\sqrt{15})}{180}   \\
           \nicefrac{1}{2}  &  \nicefrac{(10+3\sqrt{15})}{45} & \nicefrac{2}{9} &  \nicefrac{(10-3\sqrt{15})}{45}  \\
           \nicefrac{(5+\sqrt{15})}{10} & \nicefrac{(25+6\sqrt{15})}{180} & \nicefrac{(10+3\sqrt{15})}{45} & \nicefrac{5}{36}   \\
            \hline
             & \nicefrac{5}{18} & \nicefrac{4}{9}  & \nicefrac{5}{18}
        \end{array} $}
        \\
        
        \end{tabular}

    \end{table}
    \end{example}

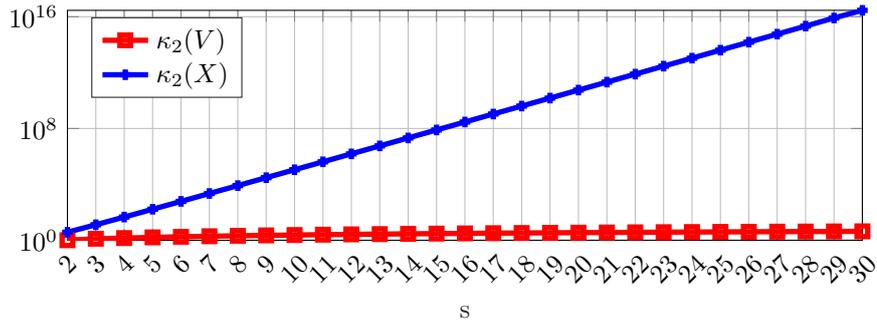
\begin{figure}[htbp]
        \centering
        \input{evcondition}
        \caption{We use the QZ algorithm to compute the left eigenvector basis for $S$, {$SV = VD$}, and $A$, {$AX = X\Lambda $}, for an increasing number of stages $s=2,\ldots,30$, and then measure their condition number.}
        \label{fig:conditioning_improvement}
    \end{figure}

    We recall that a matrix $S$ is called centroskew~\cite{CENTROSKEW} (oftentimes skew-centrosymmetric)~if
    \[
     J S J = - S,
    \]
    thus the matrix $A$ of a symmetric RK scheme is a rank-one correction of a centroskew matrix, i.e., $S = A -  \nicefrac{\mathbf{1}_s\mathbf{b}^\top}{2}$ with $J S J = - S$ and $A = S + \nicefrac{\mathbf{1}_s\mathbf{b}^\top}{2}$. This choice causes an improvement in the conditioning of the eigenvector basis as $s$ increases; see Fig.~\ref{fig:conditioning_improvement} in which we compare the conditioning of the diagonalization for the $A$ of the IRK Gauss scheme and the one of the related $S$ varying $s$.

To make use of this property, we start again from the matrix equation~\eqref{eq:rk1}, and separate its solution into two steps, one for solving the system with a centroskew matrix and one for correcting it with a low-rank matrix, namely
    \begin{eqnarray}
        \begin{array}{l|r}
        M K^{(n)} + h L K^{(n)} A^\top  =  - L \mathbf{y}^{(n)} \mathbf{1}_s^\top + F^{(n)} &  - \nonumber 
        \end{array} & \\[-0.5em]
        \begin{array}{l|r}
         M \hat{K}^{(n)} + h L \hat{K}^{(n)} S^\top  =  - L \mathbf{y}^{(n)} \mathbf{1}_s^\top + F^{(n)}\label{eq:rk1-diag} &  =\\
        \hline
        \end{array} \\[-0.35em]
        \begin{array}{l|r}
        \vspace{3mm}
        M E^{(n)} + h L E^{(n)} A^\top  =  - \frac{h}{2} [ L \hat{K}^{(n)} \mathbf{b} ] \mathbf{1}_s^\top \hspace{1.1em} & \phantom{+} \label{eq:rk1-rank}
        \end{array} 
    \end{eqnarray}
    in which $E^{(n)} = K^{(n)} - \hat{K}^{(n)}$, and~\eqref{eq:rk1-diag} can be solved through the diagonalization and parallelization procedure~\eqref{eq:diag-parallel-solve}, while~\eqref{eq:rk1-rank} is a matrix equation with a known term of rank 1 and for which, usually, the size of $L$ and $M$ is much larger than the size of $A$. We can employ a single-sided Krylov projection method with Galerkin condition~\cite{SimonciniReview} to solve this equation; see the pseudocode given in Algorithm~\ref{alg:oneside_sylvester}.

\begin{algorithm}[htb]
\KwData{Matrices $Z \in \mathbb{R}^{N \times N}$, $R \in \mathbb{R}^{s \times s}$, vectors $\mathbf{u} \in \mathbb{R}^{N}$, $\mathbf{v} \in \mathbb{R}^s$, tolerance $\epsilon$, maximum iterations $k_{\max}$}
\KwResult{Approximate solution $E$ to $Z E + E R = \mathbf{u} \mathbf{v}^\top$}

\BlankLine
Initialize $V_1 = \frac{\mathbf{u}}{\|\mathbf{u}\|}$, $\beta = \|\mathbf{u}\|$\;
\For{$j = 1, 2, \dots, k_{\max}-1$}{
    \tcc{Krylov subspace generation}
    Compute $\mathbf{w} = Z \mathbf{v}_j$\;
    \For{$i = 1, 2, \dots, j$}{
        Compute $h_{i,j} = \mathbf{v}_i^\top \mathbf{w}$\;
        Orthogonalize: $\mathbf{w} = \mathbf{w} - h_{i,j} \mathbf{v}_i$\;
    }
    Compute $h_{j+1,j} = \|\mathbf{w}\|$\;
    \If{$h_{j+1,j} = 0$}{
        \textbf{break}\;
    }
    Normalize: $\mathbf{v}_{j+1} = \frac{\mathbf{w}}{h_{j+1,j}}$\;
    Update $V_{j+1} = [\mathbf{v}_1| \cdots | \mathbf{v}_{j+1}]$\;
Form the matrix $H_j = V_j^\top Z V_j$\;
\tcc{Projected problem solution}
Solve the projected Sylvester equation $H_j Y_j + Y_j R = \beta e_1 \mathbf{v}^\top$\;
Compute the residual: $\rho_j = |h_{j+1,j}| \| \mathbf{e}_j^\top Y_j \|_F$\;
\If{$ \rho_j < \epsilon \times \| \mathbf{u} \|_F \| \mathbf{v} \|_F$}{
    \textbf{break}\;
}
}
Compute the approximate solution $E_j = V_j Y_j$\;
\caption{One-sided Krylov projection method for the matrix equation $Z E + E R = \mathbf{u} \mathbf{v}^\top$ with Arnoldi iteration}\label{alg:oneside_sylvester}
\end{algorithm}

Note in particular that the expansion of the projection space in Algorithm~\ref{alg:oneside_sylvester} requires only a matrix-vector product. In the case of a semi-discretized PDE, this operation reduces to applying the mass matrix discretization followed by solving a system involving the spatial operator. Specifically, this corresponds to computing $Z = h^{-1}L^{-1}M$. A more robust variant of Algorithm~\ref{alg:oneside_sylvester} replaces the polynomial Krylov subspace $\mathcal{K}_k(Z,\mathbf{u}) = \operatorname{span}\{ \mathbf{u},\, Z\mathbf{u},\, \ldots,\, Z^{k-1}\mathbf{u} \}$, constructed with respect to the matrix $Z$ and the vector $\mathbf{u}$, with a rational-type Krylov subspace~\cite{MR3361437}. While it is often possible to design projection spaces where the rational components have poles optimized for the problem at hand, a reliable and effective compromise between usability and efficiency is offered by the extended Krylov subspace~\cite{MR2318706}, defined as $\mathcal{E}_{2k}(Z,\mathbf{u}) = \operatorname{span}\{ \mathbf{u},\, Z\mathbf{u},\, Z^{-1}\mathbf{u},\, Z^{2}\mathbf{u},\, Z^{-2}\mathbf{u},\, \ldots,\, Z^{k-1}\mathbf{u},\, Z^{-k+1}\mathbf{u} \}$. This space captures information from both $Z$ and $Z^{-1}$, and is expanded by alternating their application to the current basis vectors. In this case, the cost of expansion increases due to the need to solve systems involving the matrix $Z$. In our setting, this corresponds to solving linear systems with the mass matrix $M$, which is typically less computationally demanding than solving systems with the spatial operator $L$.
{\begin{remark}
In the case where $L$ is singular, the solution procedure becomes more involved. 
It requires working with the Krylov subspace generated by the matrix pencil $(M,hL)$ 
and, at the projection level, employing a specialized dense solver; 
see, e.g.,~\cite{MR1383186} for a discussion of the corresponding LAPACK routines.
\end{remark}}
The computational cost of the described procedure is expressed as $s \times c_L + c_S$, where $c_L$ represents the cost of solving the linear systems involved in the block diagonalization of~\eqref{eq:rk1-diag}, and $c_S$ corresponds to the cost of solving the Sylvester equation in~\eqref{eq:rk1-rank}. By employing $n_p$ processes to simultaneously solve the step in~\eqref{eq:rk1-diag}, the cost is reduced to $\frac{s}{n_p} \times c_L + c_S$, with an additional communication overhead arising from parallel loop reduction. 

Comparing this approach with methods that rely on constructing a block-diagonalizable preconditioner for solving the entire system in~\eqref{eq:block_system} within a Krylov method, we observe a key advantage: this approach avoids the need to construct a basis for a space of size $s \times N$. Instead, it operates directly on matrices and spaces of size $N$, thus reducing memory and computational complexity.

\section{Collocation methods}\label{sec:collocation_methods}

In RK collocation methods, the solution over a single step is approximated by a polynomial, and the coefficients of this polynomial are determined by enforcing that the differential equation holds at specific points within the interval, called collocation points. The accuracy and stability of the method depend on the choice and distribution of these points. A special subset of these methods is the Gauss, Radau, and Lobatto collocation methods, which correspond to the Gauss--Legendre---the method introduced in Example~\ref{example:obtaining-coeffs}---Radau, and Lobatto quadrature points, respectively. This corresponds to the selection of the $\mathbf{b}$ and $\mathbf{c}$ as the weights and nodes of the Gauss--Legendre, Gauss--Radau, and Gauss--Lobatto quadrature formulas, respectively. For all these methods we can construct a representation in which the Butcher tableau matrix $A$, built as in~\eqref{eq:gauss2s-A-matrix}, is written as a normal matrix plus a low-rank correction, so that we can exploit the idea we explored for the case of symmetric RK methods. In the following subsection we build such a decomposition by means of the $\mathcal{W}$ transform.

\subsection{The \texorpdfstring{$\mathcal{W}$}{W} transformations}

Let $\{P_\ell\}_{\ell \geq 0}$ be the sequence of scaled and shifted degree $\ell$ Legendre polynomials on the interval $[0,1]$~\cite[P.~27]{Gautschi}, i.e., the sequence of orthonormal polynomials with respect to the $\omega(x) = 1$ measure on the $[0,1]$ interval
\[
\int_{0}^{1} P_p(x) P_q(x)\,\omega(x)\,\mathrm{d}x = \delta_{p,q} = \begin{cases}
    1, & p = q,\\
    0, & p \neq q,
\end{cases}
\]
and let $\{c_i\}_{i=1}^{s}$ be the nodes of the related Gauss-Legendre formula. We define the $\mathcal{W}_s \in \mathbb{R}^{s \times s}$ matrix
    \begin{equation}\label{eq:W-transformation-matrix}
        \left(\mathcal{W}_s\right)_{i,j} = w_{i,j} = P_{j-1}(c_i), \quad i=1,\ldots,s, \qquad j=1,\ldots,s.
    \end{equation}
    Then the following result for the Gauss method of order $2s$ holds true.
    \begin{theorem}[{\cite[Theorem~5.6]{BibleVolII}}]\label{thm:BibleVolII}
    Let $A$ be the coefficient matrix for the Gauss method of order $2s$ from~\eqref{eq:gauss2s-A-matrix} and $\mathcal{W}_s$ the matrix in~\eqref{eq:W-transformation-matrix}, then
    \[
    \mathcal{W}_s^\top A \mathcal{W}_s = \begin{bmatrix}
        \frac{1}{2} & - \xi_1 \\
        \xi_1 & 0 & -\xi_2 \\
        & \xi_2 & \ddots & \ddots \\
        & & \ddots & 0 & - \xi_{s-1} \\
        & & & \xi_{s-1} & 0 
    \end{bmatrix} = X_s, \quad \xi_k = \frac{1}{2\sqrt{4k^2 - 1}}.
    \]
    \end{theorem}
    From Theorem \ref{thm:BibleVolII} it follows that we can use the matrix $\mathcal{W}_s$ to transform problem~\eqref{eq:rk1} into a form in which we have in place of the matrix $A$ the matrix $X_s$ that is written as the sum of an antisymmetric matrix---and therefore normal---and a rank 1 correction.
    While this is not a significant addition in the case of the Gauss method, for which we can use the construction discussed in Section~\ref{sec:symmetric-RK}, this allows us to obtain effective writing for our diagonalization plus correction solution approach in other cases as well. 
    
    Consider an IRK method with tableau~\eqref{eq:butcher_tableau}, and let us define  $\mathcal{B} = \operatorname{diag}(\mathbf{b})${.} Then the matrix~\cite[pp.77--84]{BibleVolII}    \begin{equation}\label{eq:general_IRK_form}
    \begin{split}
     X_s =  &\; \mathcal{W}_s^\top \mathcal{B} A \mathcal{W}_s = \begin{bmatrix}
        \frac{1}{2} & - \xi_1 \\
        \xi_1 & 0 & -\xi_2 \\
        & \xi_2 & \ddots & \ddots \\
        & & \ddots & 0 & - \xi_{s-2} \\
        & & & \xi_{s-2} & 0 & \zeta_{s-1,s} \\
        & & & & \zeta_{s,s-1} & \zeta_{s,s}
    \end{bmatrix}, \quad \xi_k = \frac{1}{2\sqrt{4k^2 - 1}}.
    \end{split}
    \end{equation}
    The three coefficients $\zeta_{s-1,s}$, $\zeta_{s,s-1}$ and $\zeta_{s,s}$ are given in Table~\ref{tab:irk_rank} for different collocation methods.
    Moreover, we define
    \begin{equation}\label{eq:general_IRK_form_part2}
            \mathcal{D} =  \mathcal{W}_s^\top \mathcal{B} \mathcal{W}_s = \operatorname{diag}(1,1,\ldots,1,d_s),
    \end{equation}
    for which the $d_s$ is also given in Table~\ref{tab:irk_rank}.    
    \begin{table}[htbp]
        \centering
        \begin{tabular}{l|lllll}
        \toprule
             Method & $\zeta_{s,s-1}$ & $\zeta_{s-1,s}$ & $\zeta_{s,s}$ & $d_s$ & Rank \\
        \midrule
        Gauss     & $\xi_{s-1}$ & $ - \xi_{s-1}$ & $0$ & $1$ & 1 \\
        Radau IA  & $\xi_{s-1}$ & $ - \xi_{s-1}$ & $\nicefrac{1}{4s - 2}$ & $1$ & 2 \\
        Radau IIA & $\xi_{s-1}$ & $ - \xi_{s-1}$ & $\nicefrac{1}{4s - 2}$ & $1$ & 2 \\
        Lobatto IIIA & $\frac{2s-1}{s-1}\xi_{s-1}$ & 0 & 0 & $\frac{2s-1}{s-1}$ & 2 \\
        Lobatto IIIB & $0$ & $-\frac{2s-1}{s-1}\xi_{s-1}$ & 0 & $\frac{2s-1}{s-1}$ & 2 \\
        Lobatto IIIC & $\frac{2s-1}{s-1}\xi_{s-1}$ & $-\frac{2s-1}{s-1}\xi_{s-1}$ & $\frac{2s-1}{(2s-2)(s-1)}$ & $\frac{2s-1}{s-1}$ & 2 \\
        Lobatto IIIC* & $\frac{2s-1}{s-1}\xi_{s-1}$ & $-\frac{2s-1}{s-1}\xi_{s-1}$ & $-\frac{2s-1}{(2s-2)(s-1)}$ & $\frac{2s-1}{s-1}$ &  2 \\
        Lobatto IIID & $\frac{2s-1}{s-1}\xi_{s-1}$ & $-\frac{2s-1}{s-1}\xi_{s-1}$ & 0 & $\frac{2s-1}{s-1}$ & 1 \\
        \bottomrule
        \end{tabular}
                \caption{Coefficients of the $\mathcal{W}_s$-transformed Butcher tableau~\eqref{eq:general_IRK_form}-\eqref{eq:general_IRK_form_part2} for different IRK methods. The last column gives the rank of the correction needed to write $X_s$ as a skewsymmetric and hence normal matrix. Whenever $d_s \neq 1$, an additional rank 1 correction is needed to rewrite $\mathcal{D}$ as the identity.}
        \label{tab:irk_rank}
    \end{table}
    
    Let us do the change of variables $Z^{(n)} = K^{(n)} \mathcal{W}_s^{-T}$ for the equation of stages~\eqref{eq:rk1}, and multiply on the right the equation by $\mathcal{B}\mathcal{W}_s$, hence
    \[
    M Z^{(n)} \mathcal{W}_s^\top \mathcal{B}\mathcal{W}_s + h L Z^{(n)} \mathcal{W}_s^\top A^\top \mathcal{B}\mathcal{W}_s = - L \mathbf{y}^{(n)} \mathbf{b}^\top\mathcal{W}_s + F^{(n)} \mathcal{B}\mathcal{W}_s,\\
    \]
    giving the transformed generalized Sylvester matrix-equation
    \[
    \begin{cases}
    M Z^{(n)} \mathcal{D} + h L Z^{(n)} X_s^\top = - L \mathbf{y}^{(n)} \mathbf{b}^\top\mathcal{W}_s + F^{(n)} \mathcal{B}\mathcal{W}_s,\\
    K^{(n)} = Z^{(n)} \mathcal{W}_s^\top.
    \end{cases}
    \]

To apply the same method adopted for the symmetric schemes in~\eqref{eq:rk1-diag}-\eqref{eq:rk1-rank}, we need to explicitly define a low-rank transformation that diagonalizes the system, i.e., makes $X_s$ skew-symmetric, and ensures the diagonal $\mathcal{D}$ is corrected to the identity. This is achieved by setting\footnote{For all cases in Table~\ref{tab:irk_rank} except Lobatto IIIB for which we substitute $C_2 = [\mathbf{e}_1,\mathbf{e}_s]$ with  $C_2 = [\mathbf{e}_1,\mathbf{e}_{s-1}]$, observe also that for the Lobatto IIIA and IIIB $\zeta_{s,s} = 0$. This choice avoids a zero eigenvalue in $\hat{X}_s$.}:
\[
\begin{split}
\hat{X}_s &= X_s + \begin{bmatrix}
    -\nicefrac{1}{2} & 0 \\
    0 & \vdots \\
    \vdots & 0 \\
    0 & -\zeta_{s,s-1} - \zeta_{s-1,s} \\ %
    0 & - \zeta_{s,s}
\end{bmatrix} 
\begin{bmatrix}
    1 & 0 \\
    0 & \vdots \\
    \vdots & 0 \\
    0 & 1
\end{bmatrix}^\top = X_s + C_1 C_2^\top, \\
I_s &= \mathcal{D} + (1-d_{s,s}) \mathbf{e}_s \mathbf{e}_s^\top.
\end{split}
\]

We solve the block-diagonal linear system with $\hat{X}_s$ in place of $X_s$ and $I_s$ in place of~$\mathcal{D}$:
\[
\left(I_s \otimes M + h \hat{X}_s \otimes L\right) \operatorname{vec}\left(\hat{Z}^{(n)}\right) = \operatorname{vec}\left(-L \mathbf{y}^{(n)} \mathbf{b}^\top \mathcal{W}_s + F^{(n)} \mathcal{B} \mathcal{W}_s\right),
\]
where $\hat{X}_s$ is diagonalizable as $\hat{X}_s = Q \Lambda Q^H$ with $Q^H Q = I_s$. Consequently, we can rewrite:
\[
(Q \otimes I_N) \left(I_s \otimes M + h \Lambda \otimes L\right) (Q^H \otimes I_N) \operatorname{vec}\left(\hat{Z}^{(n)}\right) = \operatorname{vec}\left(-L \mathbf{y}^{(n)} \mathbf{b}^\top \mathcal{W}_s + F^{(n)} \mathcal{B} \mathcal{W}_s\right).
\]

The correction matrix equation can then be written as:
\begin{equation}\label{eq:generic-kutta-correction}
\begin{array}{rcl|l}
M {Z}^{(n)} \mathcal{D} + h L {Z}^{(n)} {X}_s^\top & = & - L \mathbf{y}^{(n)} \mathbf{b}^\top\mathcal{W}_s + F^{(n)} \mathcal{B}\mathcal{W}_s  & - \\
M \hat{Z}^{(n)} + h L \hat{Z}^{(n)} \hat{X}_s^\top & = & - L \mathbf{y}^{(n)} \mathbf{b}^\top\mathcal{W}_s + F^{(n)} \mathcal{B}\mathcal{W}_s   & = \\
\hline\\[-0.9em]
M E^{(n)} \mathcal{D} + h L E^{(n)} X_s^\top & = & (1-d_{s,s})M\hat{Z}^{(n)}\mathbf{e}_s \mathbf{e}_s^\top + hL\hat{Z}^{(n)}C_2C_1^\top
\end{array}
\end{equation}
where $E^{(n)} = Z^{(n)} - \hat{Z}^{(n)}$. The right-hand side is now a matrix of rank two or three, depending on the linear independence of the vectors forming the three dyads.

\begin{remark}
To obtain the minimal rank version of the right-hand side of~\eqref{eq:generic-kutta-correction}, it is sufficient to collect the three vectors together and apply the Gram-Schmidt procedure to the columns and rows respectively to discard the only possible linearly dependent vector. To solve the matrix equation~\eqref{eq:generic-kutta-correction} it is then possible to apply the Algorithm~\ref{alg:oneside_sylvester} in which the Krylov space is replaced by the block Krylov space, i.e., block Arnoldi~\cite[Algorithm~6.22]{SAAD}, with block size equal to the rank of the right-hand side.
\end{remark}

\begin{remark}
Note that for an odd number of stages $s$, the center-skew matrix $S$ or the skew-symmetric matrix $\hat{X}_s$ will always have a zero eigenvalue by construction; see the examples in Fig.~\ref{fig:eigentrouble}. When the mass matrix is either the identity matrix or a more general non-singular matrix, choosing an odd number of stages does not pose any issues for the solvability of the shifted systems.  
\begin{figure}[htbp]
    \centering
    \includegraphics[width=\columnwidth]{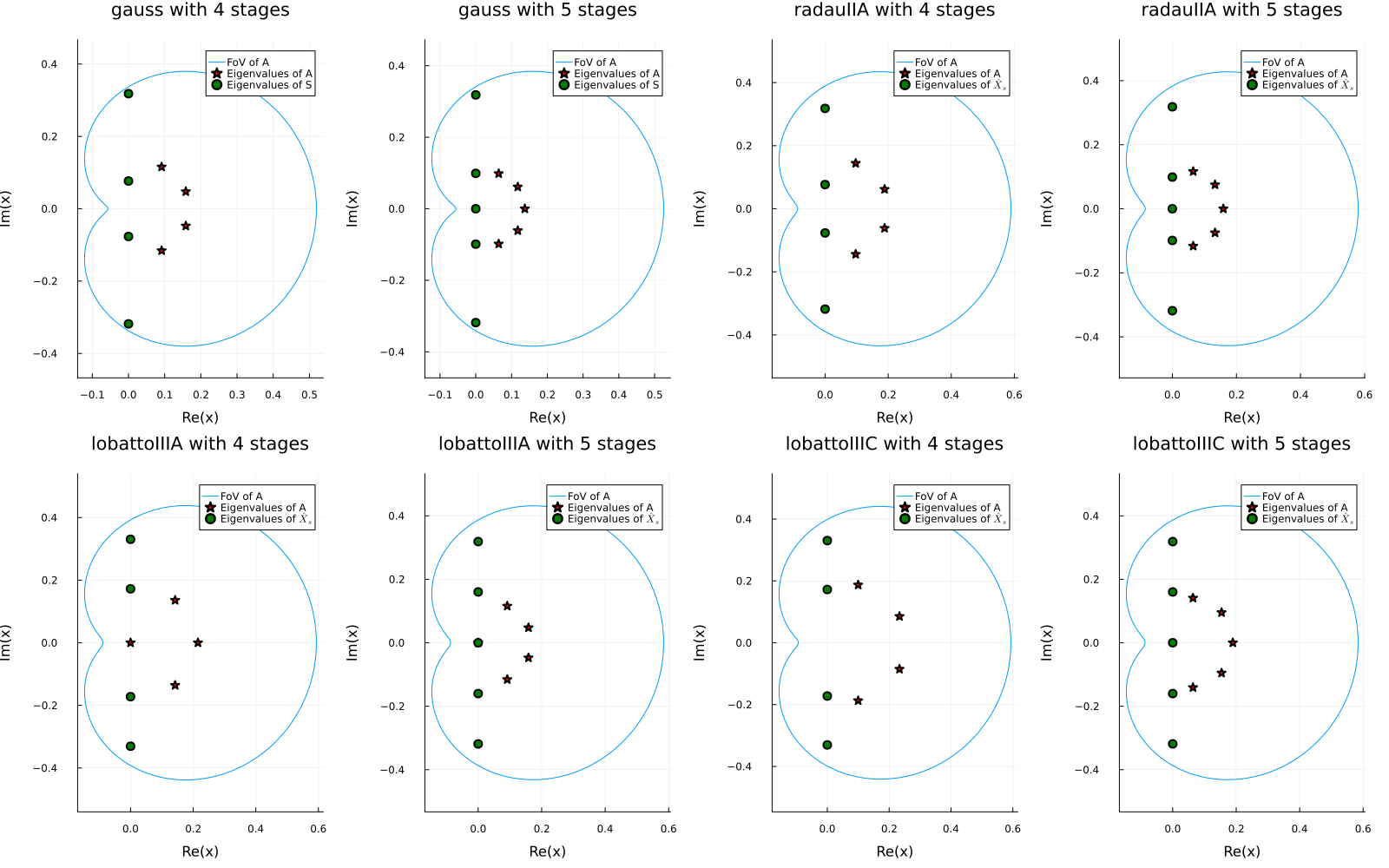}
    \caption{Visualization of the eigenvalues ($\star$) and the field of values of the Butcher matrix $A$ for various Runge--Kutta schemes, along with the eigenvalues of the perturbed matrices $S$ and $\hat{X}_s$ ($\circ$) for stage numbers $s = 4, 5$.}
    \label{fig:eigentrouble}
\end{figure}
However, if the considered case involves a differential-algebraic equation (DAE) or a singular mass matrix, one of the systems to be solved using this approach will have a non trivial kernel.

\end{remark}

\section{The nonlinear case: the simplified Newton method}\label{sec:nonlinear}

In the nonlinear case, the differential equation in~\eqref{eq:linear-pde} is substituted by
    \begin{equation}\label{eq:nonlinearsystem}
            \begin{cases}
        \frac{\partial u}{\partial t} + \varTheta[u,t] = 0, & (\mathbf{x},t) \in \Omega \times (0,T] , \quad \Omega \subseteq \mathbb{R}^{d},   \\
    u(0) = u_0(\mathbf{x}) , \\
    \mathcal{B}\left[ u(\mathbf{x},t) \right] = 0, &  (\mathbf{x},t) \in \partial \Omega \times [0,T],
    \end{cases},\quad d \in \mathbb{N},
    \end{equation}
    for $\varTheta[\cdot,\cdot]$ a---possibly non autonomous---nonlinear operator involving partial derivatives of the function $u$. After a suitable space discretization, we come back to the case 
    \[
    \begin{cases}
        M\mathbf{u}{'}(t) = - \Theta(\mathbf{u},t), & t \in (0,T] ,   \\
    \mathbf{u}(0) = \mathbf{u},
    \end{cases} \qquad M \in \mathbb{R}^{N \times N}, \; \Theta : \mathbb{R}^N \times [0,T] \rightarrow \mathbb{R}^N.
    \]
    To solve the stage equations in~\eqref{eq:stages_Runge_Kutta} we need to apply a suitable linearization procedure to it, typically a Newton-like or a Picard iteration. In the Newton case, we start by introducing the auxiliary functions
    \begin{equation}\label{eq:hat-theta}
    \hat{\Theta}_i(\mathbf{k}_i^{(n)}) \equiv M \mathbf{k}_{i}^{(n)} - M \mathbf{u}^{(n)} + h \sum_{j=1}^{s} a_{i,j} \Theta\left( \mathbf{k}_{j}^{(n)}, t_n + c_j h\right) = 0, \qquad i=1,\ldots,s,
    \end{equation}
    and thus the vector function $\hat{\boldsymbol{\Theta}}(\mathbf{k},t) = [ \hat{\Theta}_1,\ldots,\hat{\Theta}_s ]^\top$. From the solution of~\eqref{eq:hat-theta} we obtain the expression for the stages with which we advance of a time-step by
    \begin{equation}\label{eq:time_advance_nonlinear}
    M \mathbf{u}^{(n+1)} = M \mathbf{u}^{(n)} - h \sum_{i=1}^{s} b_i \Theta( \mathbf{k}_j^{(n)}, t_n + c_i h).      
    \end{equation}
    The full Newton iteration for the approximation of the zeros of~\eqref{eq:hat-theta} is given by
    \begin{equation}\label{eq:full_newton}
                \begin{cases}
            \text{Given }\boldsymbol{\kappa}^0 = [\mathbf{k}_1^{(n),0},\ldots,\mathbf{k}_s^{(n),0}],\\
            J_{\hat{\boldsymbol{\Theta}}}( \boldsymbol{\kappa}^p ) \mathbf{d}^p = \hat{\boldsymbol{\Theta}}(\boldsymbol{\kappa}^p), & p \geq 0\\
            \boldsymbol{\kappa}^{p+1} = \boldsymbol{\kappa}^{p} - \mathbf{d}^p,%
        \end{cases}
    \end{equation}
    for $J_{\hat{\boldsymbol{\Theta}}}( \boldsymbol{\kappa}^p )$ a suitable approximation of the Jacobian of the vector function $\hat{\boldsymbol{\Theta}}$ at the current step, and which can be expressed as
    \[
    \begin{split}
    J_{\hat{\boldsymbol{\Theta}}}( \boldsymbol{\kappa}^p ) %
    = & \left(  A^{-1} \otimes M + h \begin{bmatrix}
    J_{\Theta}^{(1)} \\
         & \ddots \\
     & & J_{\Theta}^{(s)}
    \end{bmatrix} \right)( A \otimes I_N )
    \end{split}
    \]
    where $J_{\Theta}^{(i)}$, $i=1,\ldots,s$, is the Jacobian of the $\Theta$ function with respect to the $\mathbf{k}_i^{(n)}$ variables. %
    To be able to return to the linear case, we must assume that all Jacobian matrices are computed for the same time-step, that is, that the system is autonomous, or close to being so; i.e., we build what is called a \emph{simplified Newton} method~\cite{SimplifiedNewton,MR1790038}. In such a case we build the single matrix $J_{\Theta}$ for all the stages and use
    \begin{equation*}
    J_{\hat{\boldsymbol{\Theta}}} = \left(   I_s  \otimes M + h A \otimes J_{\Theta} \right),
    \end{equation*}
    in~\eqref{eq:full_newton}; e.g., as in~\cite{MR4167091}, we can select $J_{\Theta}$ as the arithmetic average of the $J_{\Theta}^{(i)}$, $i=1,\ldots,s$. Now, to solve for the Newton direction we apply the correction procedure. In summary, to advance from time-step $t_{n}$ to time-step $t_{n+1}$ we use~\eqref{eq:time_advance_nonlinear}, where $\mathbf{k}_i^{(n)}$ is obtained from the simplified Newton method applied to $\hat{\boldsymbol{\Theta}}(\mathbf{k}_i^{(n)}) = 0$, i.e., we select $\boldsymbol{\kappa}^0 = [\mathbf{k}_1^{(n),0},\ldots,\mathbf{k}_s^{(n),0}]$ and solve
    \[
    \begin{split}
        J_{\hat{\boldsymbol{\Theta}}} \mathbf{d}^p = \hat{\boldsymbol{\Theta}}(\boldsymbol{\kappa}^p), \\
        \left(   I_s  \otimes M + h A \otimes J_{\Theta} \right) \mathbf{d}^p =  \hat{\boldsymbol{\Theta}}(\boldsymbol{\kappa}^p).
    \end{split}
    \]
    To simplify the discussion, we assume working with a symmetric scheme as in Section~\ref{sec:symmetric-RK}, and solve instead
    \[
    \left(   I_s  \otimes M + h S \otimes J_{\Theta} \right) \hat{\mathbf{d}}^p =  \hat{\boldsymbol{\Theta}}(\boldsymbol{\kappa}^p),
    \]
    for $S = A - \nicefrac{\mathbf{1}_s \mathbf{b}^\top}{2}$, by block-diagonalization
    \begin{equation*}  
        (X \otimes I_N) \left(   I_s  \otimes M + h \Lambda \otimes J_{\Theta} \right)  (X^{-1} \otimes I_N) \hat{\mathbf{d}}^p   =    \hat{\boldsymbol{\Theta}}(\boldsymbol{\kappa}^p),
    \end{equation*}
    \begin{eqnarray}
    \left(   I_s  \otimes M + h \Lambda \otimes J_{\Theta} \right)  \mathbf{z} = &\;  (X^{-1} \otimes I_N) \hat{\boldsymbol{\Theta}}(\boldsymbol{\kappa}^p), \label{eq:linear_systems_from_non_linear1} \\
    \hat{\mathbf{d}}^p = &\; ( X \otimes I_N ) \mathbf{z}. \nonumber
        \end{eqnarray}
    If we denote by $\Delta_p$, $\hat{\Delta}_p$ and $\theta$ the matricization of size $N \times s$ of the vectors $\mathbf{d}^p$, $\hat{\mathbf{d}}^{p}$ and $\hat{\boldsymbol{\Theta}}(\boldsymbol{\kappa}^p)$ respectively, then we recover the Newton direction $\mathbf{d}^p$ via the solution of the Sylvester equation
    \begin{equation}\label{eq:sylvester_for_newton}
            \begin{array}{rl|c}
    M \Delta_p + h J_{\Theta} \Delta_p A^\top = & \theta & - \\
    M \hat{\Delta}_p + h J_{\Theta} \hat{\Delta}_p S^\top = & \theta & =\\
    \midrule
    M E_p + h J_{\Theta} E_p A^\top = & - \frac{h}{2}\left(J_{\Theta} \hat{\Delta}_p \mathbf{b} \right) \mathbf{1}_s^\top
    \end{array},
    \end{equation}
    with $E_p = \Delta_p - \hat{\Delta}_p$. Hence the Newton update is computed as
    \[
    \boldsymbol{\kappa}^{p+1} = \boldsymbol{\kappa}^p - \operatorname{vec}\left( \hat{\Delta}_p + E_p \right).
    \]
    The procedure is analogous for the collocation methods discussed in Section~\ref{sec:collocation_methods}.

\section{Implementation details and numerical examples}\label{sec:code_and_implementation}
The approach discussed in this work is implemented in Julia and available at the GitHub repository \href{https://github.com/Cirdans-Home/SP\_IRK.jl}{Cirdans-Home/SP\_IRK.jl}. Our implementation takes full advantage of Julia's shared-memory parallelism to accelerate critical components of the computation. One of the central tasks is the parallel computation of LU factorizations required for solving the linear systems that appear in equations~\eqref{eq:linear_systems}, and~\eqref{eq:linear_systems_from_non_linear1}. This is accomplished by concurrently factorizing matrices of the form \(M + h \lambda_j L\)---respectively \(M + h \lambda_j J_\Theta\) and $J_\Theta$ a suitable approximation of the Jacobian. The operation is encoded as follows:
\begin{lstlisting}[language=Julia]
lhs = fetch.([Threads.@spawn factorize(Mass + h*ev[j]*L) for j in 1:s])
\end{lstlisting}
In this snippet, the LU factorization of each matrix is computed in a separate asynchronous task using \texttt{Threads.@spawn}, and the resulting factorizations are aggregated into the array \texttt{lhs} via \texttt{fetch}. For autonomous linear problems, where the matrix remains unchanged over time, this factorization is computed just once and then reused for all time-steps. This step is efficiently handled by the external UMFPACK library~\cite{MR2075981} interfaced through Julia.

In addition to the factorization, the code also exploits multi-threading for assembling the right-hand side (RHS) of the linear systems and for applying the computed factorizations to obtain the solutions. The assembly of the RHS involves evaluating the ODE vector field at modified time points, and the corresponding loop is parallelized using \texttt{Threads.@threads} as shown below for the linear case~\eqref{eq:linear_problem}:
\begin{lstlisting}[language=Julia]
Threads.@threads for j in eachindex(c)
    @inbounds F[:,j] = f(ti + c[j]*h)
end
\end{lstlisting}
Here, the use of \texttt{@inbounds} eliminates array bounds checking, thereby reducing overhead since the dimensions of \lstinline[language=Julia]{F} have already been verified. Following the assembly, the precomputed LU factorizations are applied to solve the linear systems in parallel:
\begin{lstlisting}[language=Julia]
Threads.@threads for j in eachindex(lhs)
     @inbounds K[:,j] = lhs[j] \ r[:,j]
end
\end{lstlisting}
Each thread independently applies the corresponding factorization from \texttt{lhs} to the RHS vector segment in \texttt{r}, storing the resulting solution in \texttt{K} while again employing \texttt{@inbounds} for performance optimization. It is important to note that while multi-threading significantly accelerates the computation by exploiting parallel execution, there is an inherent overhead associated with creating and managing the threaded environment. This overhead stems from the cost of thread creation, task scheduling, which may become non-negligible, particularly for smaller problem sizes or less intensive computational tasks. Therefore, when applying parallelization, it is critical to balance these overhead costs against the performance gains to ensure that the overall efficiency is maintained. It is also worth mentioning that another bottleneck for parallelism arises in the correction step following the solution of the linear systems, where the Sylvester matrix equation~\eqref{eq:rk1-rank}---respectively~\eqref{eq:generic-kutta-correction} and~\eqref{eq:sylvester_for_newton}---is solved. This step, due to its inherent sequential nature, limits further speedups that might be obtained through parallelization. Furthermore, the maximum number of threads that can be effectively utilized is constrained by the number $s$ of stages in the underlying IRK method, as each stage represents an independent computation whose total number sets the upper bound for concurrent execution.

\subsection{Numerical examples}\label{sec:numerical_examples}

We examine two differential problems to demonstrate the effective application of the strategy discussed earlier. Our focus lies on a framework that addresses PDEs using the longitudinal method of lines. This approach begins by discretizing the spatial domain, converting the PDE into a system of ODEs, which can be linear or nonlinear. These ODEs are subsequently solved using the proposed integration methods. The numerical experiments contained in this section are conducted on a single node of the Toeplitz cluster, located at the Green Data Center of the University of Pisa. The node is equipped with an AMD EPYC 7763 64-Core Processor CPU, featuring 2 threads per core, 64 cores per socket, 2 sockets, and a total of 2~\si{\tera\byte} of RAM, which makes it particularly suitable for testing the multi-threaded programming model approach discussed above.

We analyze two specific cases: a linear heat transport problem in Section~\ref{sec:heat-transport}, and a nonlinear wave problem in Section~\ref{sec:nonlinear-wave}. These examples test the methods introduced in Section~\ref{sec:collocation_methods} and \ref{sec:nonlinear}, respectively.

\subsubsection{An heat transport problem}\label{sec:heat-transport}

\begin{figure}[htbp]
    \centering

    \begin{subfigure}[t]{0.42\columnwidth}
    \includegraphics[width=\columnwidth]{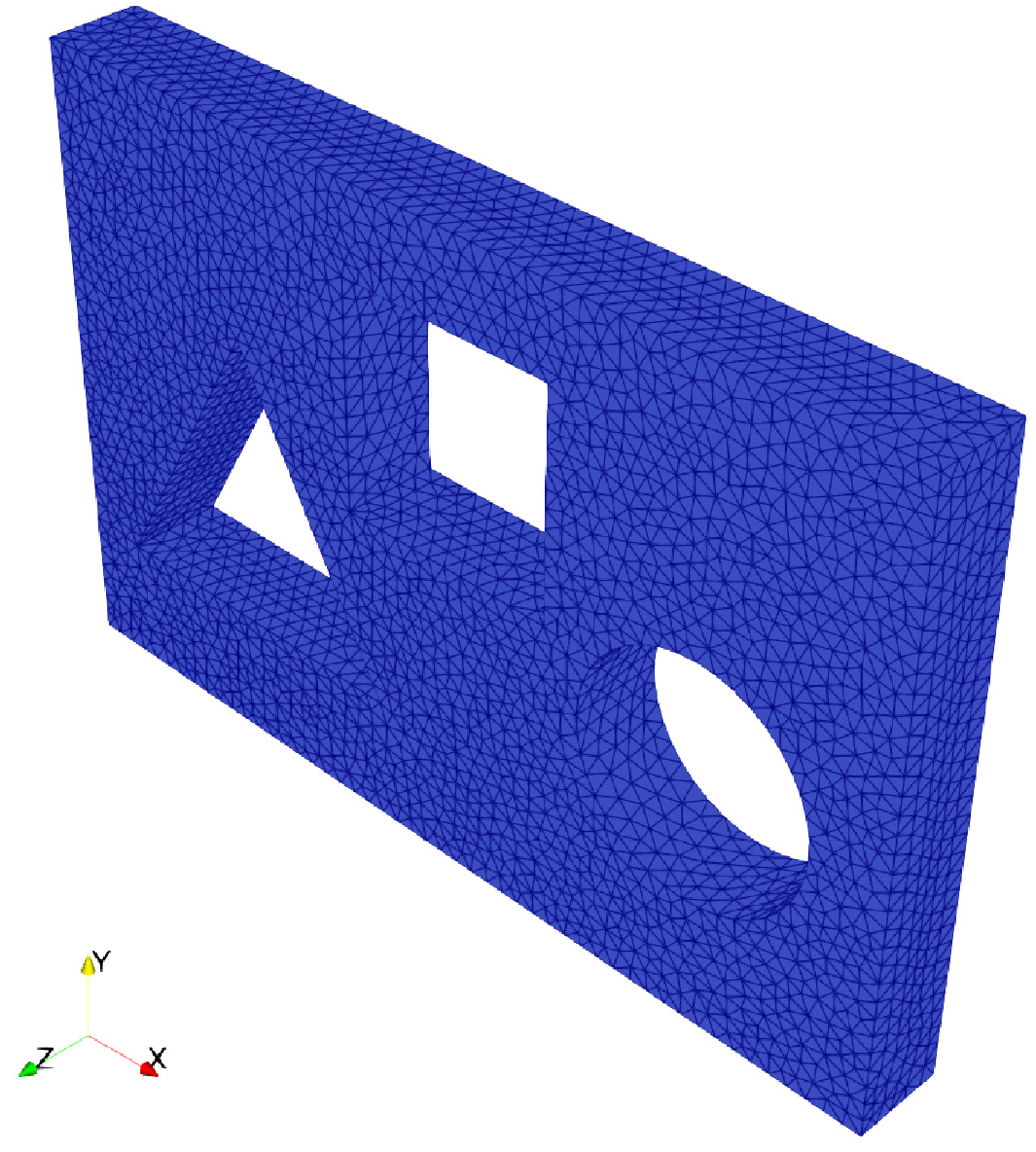}

    \caption{Domain with an example mesh for the heat transport problem~\eqref{eq:heat-transport-problem}.}\label{fig:heat_problem_domain}
    \end{subfigure}\hfil
    \begin{subfigure}[t]{0.42\columnwidth}
    \includegraphics[width=\columnwidth]{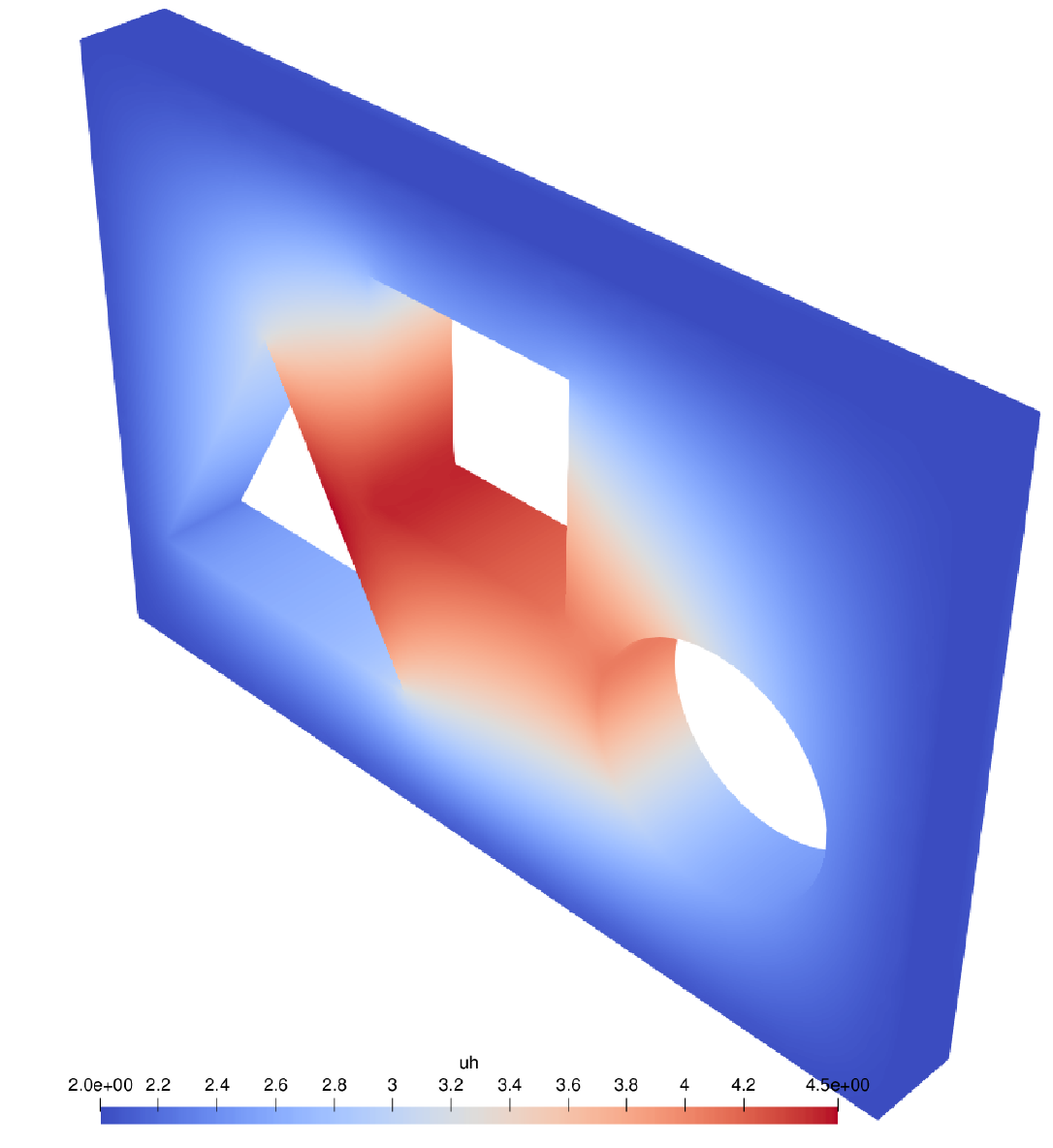}

    \caption{Solution at $T = 1$ for the choice of parameters in~\eqref{eq:heat-problem-coefficients}.\label{fig:heat_prob_sol}}
    \end{subfigure}
    
    \caption{Heat transport problem from Section~\ref{sec:heat-transport}.}
    \label{fig:solution_heat_transport}
\end{figure}

 We solve a finite element discretization of the heat equation on the 3D domain depicted in Fig.~\ref{fig:solution_heat_transport} with mixed Dirichlet and Neumann boundary conditions. Dirichlet boundary conditions are applied on the outer sides of the holed prism $\Gamma_{\mathfrak{D}}$, while non-homogeneous Neumann conditions are applied to the three internal boundaries of the holes ($\Gamma_\alpha$ for the triangular hole, $\Gamma_\gamma$ for the square hole, and $\Gamma_\beta$ for the circular hole). Homogeneous Neumann boundary conditions are applied in the remaining portion of the boundary $\Gamma_W$---hence $\Gamma_{\mathfrak{N}} = \Gamma_\alpha \cup \Gamma_\gamma \cup \Gamma_\beta \cup \Gamma_W$:
 \begin{equation}\label{eq:heat-transport-problem}
\begin{cases}
    \frac{\partial u}{\partial t} - \nabla^2 u = f(\mathbf{x},t), & \text{ in } \Omega\times(0,T], \\
    u = g_1(\mathbf{x}), & \text{ on } \Gamma_{\mathfrak{D}}\times(0,T],\\
    \nabla u \cdot \vec{\mathbf{n}} = g_2(\mathbf{x}), & \text{ on } \Gamma_{\mathfrak{N}}\times(0,T],\\
    u(\mathbf{x},0) = u_0(\mathbf{x}),
\end{cases}
 \end{equation}
 where $\vec{\mathbf{n}}$ is the outward pointing normal with respect to the Neumann boundaries; see Fig.~\ref{fig:heat_problem_domain} for a depiction of the domain.

To arrive at the semi-discrete formulation in~\eqref{eq:linear_problem} we employ the Finite Element \texttt{Gridap.jl} library~\cite{Badia2020,Verdugo2022} and select
\begin{equation}\label{eq:heat-problem-coefficients}
    f(\mathbf{x},t) = (1.5-x_1)(1-x_2)(1-x_3) + \frac{1}{2}\sin(2\pi\,t), \quad g_1(\mathbf{x}) = 2, \quad g_2(\mathbf{x}) = 3.
\end{equation}
For all the following examples we set $T = 1$, and $h = 10^{-2}$; see Fig.~\ref{fig:heat_prob_sol} for a depiction of the solution.

\begin{figure}[htb]
    \centering
    \begin{subfigure}{0.45\columnwidth}
    \includegraphics[width=\columnwidth]{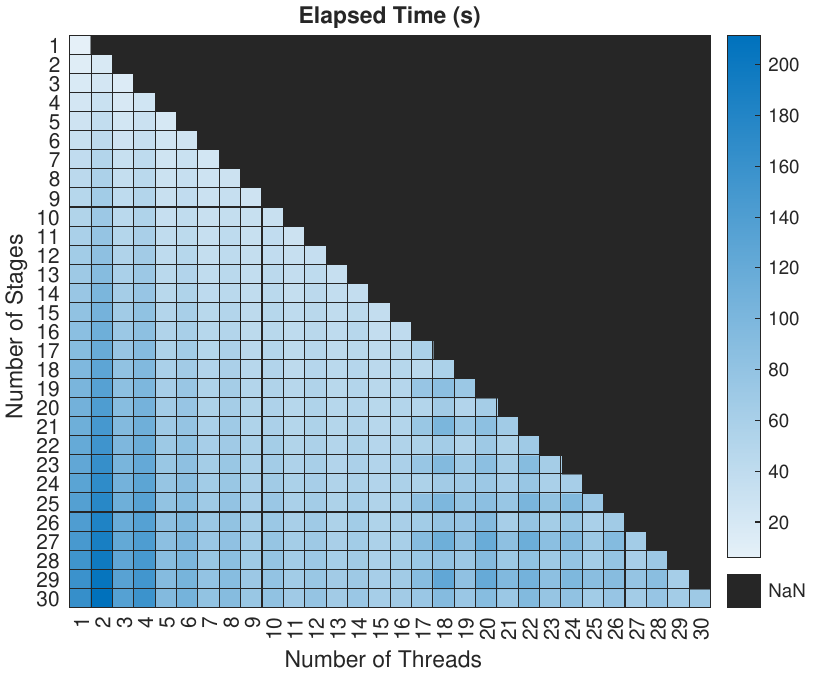}
    
    \caption{Elapsed time (s)\label{fig:elapsed_time}}
    \end{subfigure}
    \begin{subfigure}{0.45\columnwidth}
    \includegraphics[width=\columnwidth]{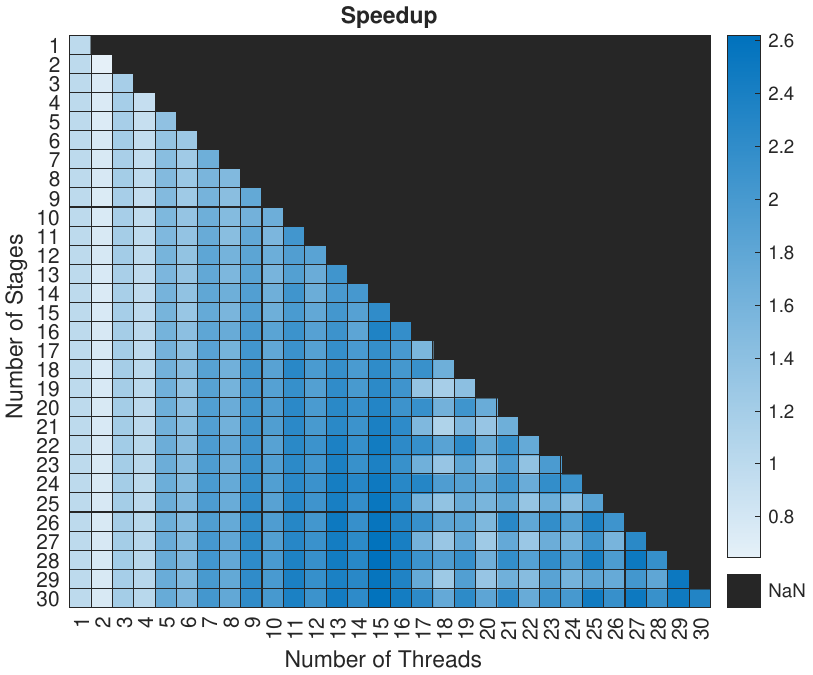}

    \caption{Speedup \label{fig:speedup}}
    \end{subfigure}
    
    \caption{Elapsed time and speedup for the heat transport problem in Section~\ref{sec:heat-transport} obtained by employing the stage-parallel procedure with the Gauss symmetric scheme and number of stages and threads running from $1$ to $30$.}
    \label{fig:heat_gauss}
\end{figure}
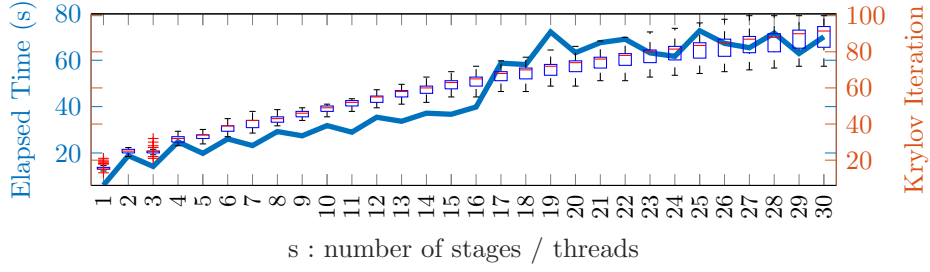
\begin{figure}[htb]
    \centering
    \input{iteration_time}
    
    \caption{Elapsed time and iteration scaling with respect to the number of stages and equal number of threads for the solution of the Sylvester equation~\eqref{eq:rk1-rank} in the case of the Gauss scheme, implemented as a symmetric scheme.}
    \label{fig:time_and_krylov}
\end{figure}

Fig.~\ref{fig:elapsed_time} illustrates the elapsed time for the solution of the Gauss IRK method using a stage-parallel approach. As the number of threads increases, a noticeable reduction in elapsed time is observed, highlighting the effectiveness of parallelizing across stages. However, the rate of improvement diminishes beyond a certain point, suggesting that the overhead caused by the handling of the multi-thread environment in Julia and the solution of the associated matrix equation starts to dominate at higher thread counts. This behavior is further quantified in Fig.~\ref{fig:speedup}, which shows the parallel speedup achieved. The speedup increases significantly to around $3\times$ with the number of threads, but again, exhibits sublinear scaling as the number of computing units grows. To better investigate the numerical behavior of the procedure, in Fig.~\ref{fig:time_and_krylov} we focus on the solve time for the number of stages equal to the number of threads, i.e., to the diagonal of Fig.~\ref{fig:elapsed_time}. The left axis reports the total elapsed time measured in seconds, while the right axis shows the Krylov iteration count. We observe that the number of iterations performed by the polynomial Krylov method described in Algorithm~\ref{alg:oneside_sylvester} grows with $s$. Indeed, if we interpret the scaling as a weak scaling with $s$, we observe that increasing the problem size from $2N$ to $30N$ shows $3.7\times$ increase in computational time. Overall, the results confirm that the stage-parallel strategy provides a practical mean to enhance computational efficiency for time integration in stiff systems.

\begin{figure}[htb]
    \centering
    \begin{subfigure}{0.45\columnwidth}
    \includegraphics[width=\columnwidth]{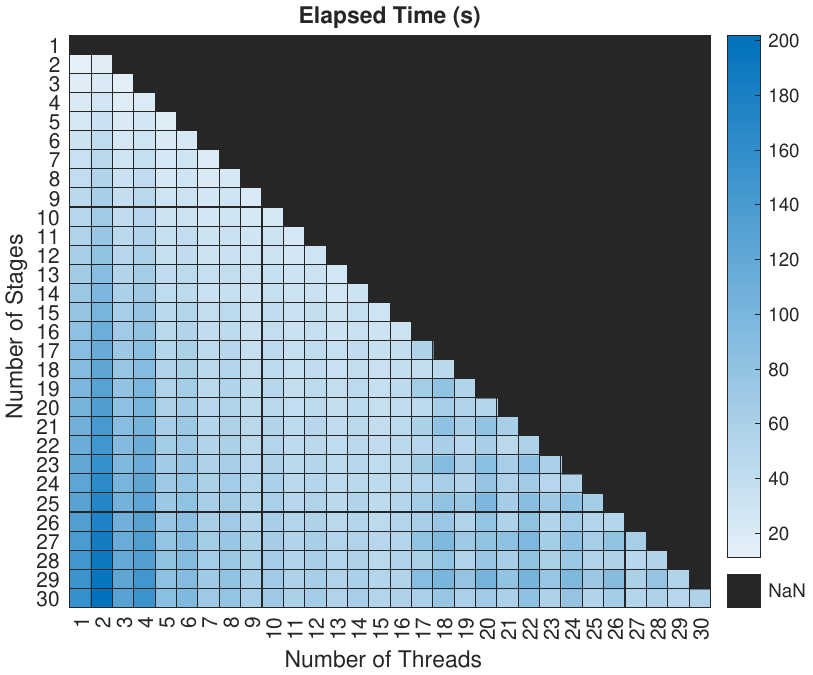}
    
    \caption{Elapsed time (s)\label{fig:elapsed_time_radau}}
    \end{subfigure}
    \begin{subfigure}{0.45\columnwidth}
    \includegraphics[width=\columnwidth]{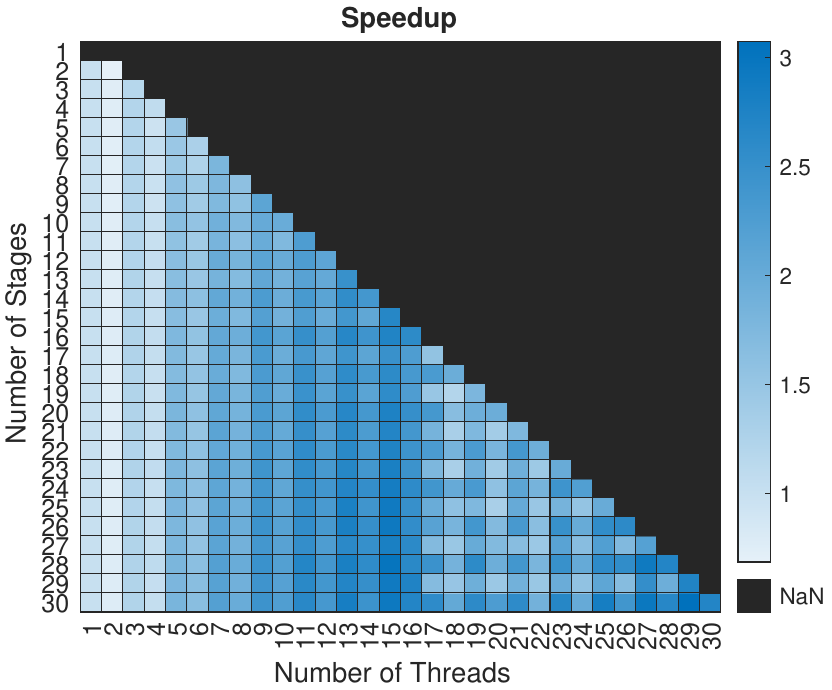}

    \caption{Speedup \label{fig:speedup_radau}}
    \end{subfigure}
    
    \caption{Elapsed time and speedup for the heat transport problem in Section~\ref{sec:heat-transport} obtained by employing the stage-parallel procedure with the Radau IIA collocation scheme and number of stages  running from $2$ to $30$ and threads running from $1$ to $30$.}
    \label{fig:heat_radauIIA}
\end{figure}

\begin{figure}[htb]
    \centering
    \input{iteration_time_radau}
    
    \caption{Elapsed time and iteration scaling with respect to the number of stages and equal number of threads for the solution of the Sylvester equation~\eqref{eq:generic-kutta-correction} in the case of the Radau IIA scheme, implemented as a collocation scheme.}
    \label{fig:time_and_krylov_radau}
\end{figure}
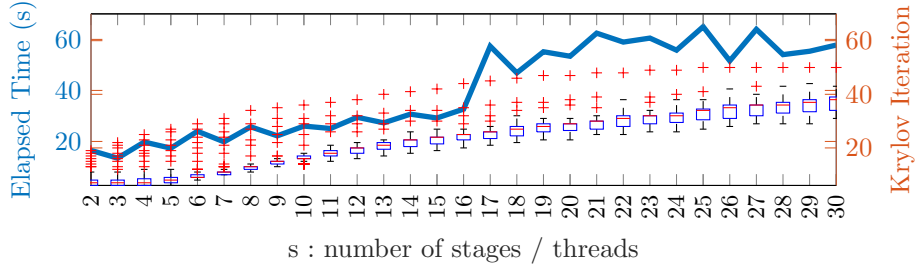

We repeat the experiment using the same configuration as before, but replacing the Gauss method with the Radau IIA method, which is based on its interpretation as a collocation method (Section~\ref{sec:collocation_methods}).  The elapsed time and computed speedup are presented in Fig.~\ref{fig:heat_radauIIA}, and can be compared with the results shown in Fig.~\ref{fig:heat_gauss}.
The behavior in both cases is analogous, with the Radau IIA scheme\footnote{Note that for Radau IIA the scheme with $s=1$ is meaningless.} achieving a slightly higher maximum speedup compared to the Gauss method. This improvement can be attributed to the more efficient solution of the Sylvester equation~\eqref{eq:generic-kutta-correction} in the Radau IIA case, as opposed to the Gauss case where the alternative form~\eqref{eq:rk1-rank} was used instead. See Fig.~\ref{fig:time_and_krylov_radau} for the Radau IIA results, and compare them with the corresponding Fig.~\ref{fig:time_and_krylov} for the Gauss case.

\subsubsection{A nonlinear wave equation}\label{sec:nonlinear-wave}

\begin{figure}[htbp]
    \centering

    \begin{subfigure}{0.45\columnwidth}
    \includegraphics[width=\columnwidth]{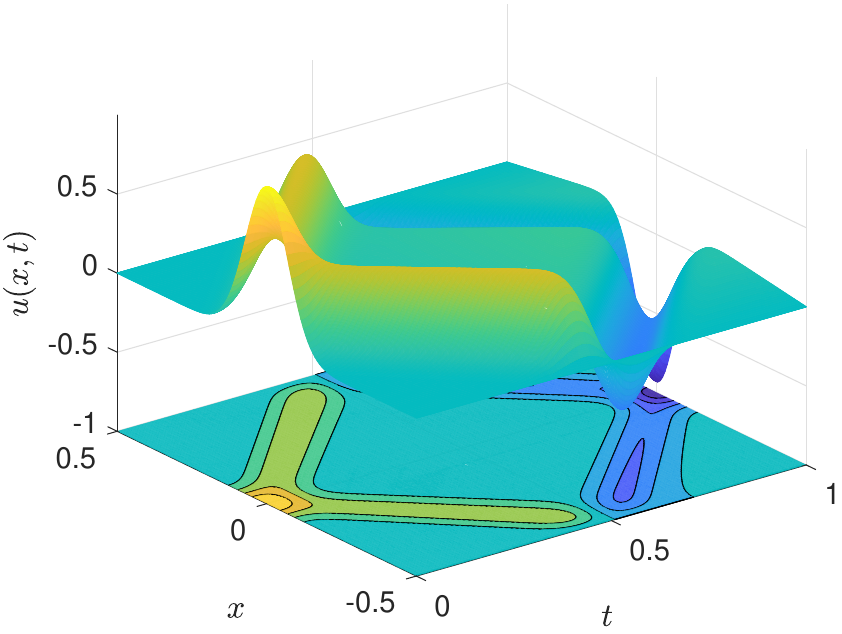}
    
    \caption{$u(x,t)$, $x \in [\nicefrac{1}{2},\nicefrac{1}{2}]$, $t \in [0,1]$}
    \end{subfigure}
    \begin{subfigure}{0.45\columnwidth}
    \includegraphics[width=\columnwidth]{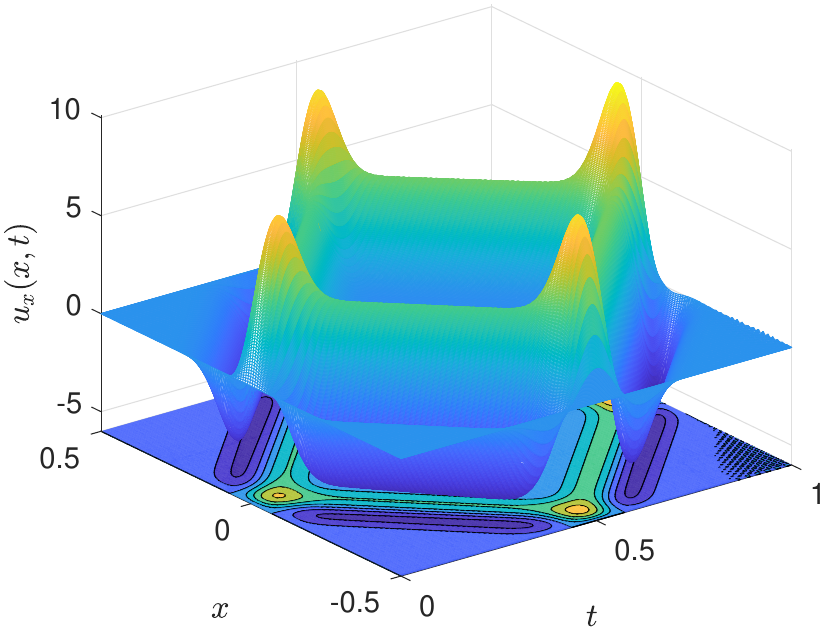}
    
    \caption{$u_x(x,t)$, $x \in [\nicefrac{1}{2},\nicefrac{1}{2}]$, $t \in [0,1]$}
    \end{subfigure}
    
    \caption{Solution of the nonlinear wave equation~\eqref{eq:nonlinearwaeve} with $\beta = 10$, $u(x,0) = \exp(-100 x^2)$, $u_x(x,0) = 0$ obtained with $\Delta x = \nicefrac{1}{127}$ and $h = \nicefrac{1}{1270}$ using the symmetric Gauss scheme with $s = 4$ stages.}
    \label{fig:solution-of-wave}
\end{figure}

We next consider the nonlinear wave equation from~\cite[\S 5.2]{MR4167091},
\begin{equation}\label{eq:nonlinearwaeve}
\left\{
\begin{aligned}
\partial_{tt} u &= \partial_{xx} u + \beta u^2, & (x,t) &\in \left(-\tfrac{1}{2}, \tfrac{1}{2}\right) \times (0,1), \\
u\left(-\tfrac{1}{2}, t\right) &= u\left(\tfrac{1}{2}, t\right) = 0, & t &\in (0,1), \\
u(x,0) &= e^{-100x^2},\quad u_t(x,0) = 0, & x &\in \left(-\tfrac{1}{2}, \tfrac{1}{2}\right).
\end{aligned}
\right.
\end{equation}
The parameter $\beta > 0$ represents the strength of the nonlinearity. We semi-discretize in space the equation by employing centered finite differences as in~\cite{MR4167091} with a mesh-size $\Delta x$, which results in the system of ODEs of the form~\eqref{eq:nonlinearsystem}
\[
\begin{bmatrix}
\mathbf{u} \\
\mathbf{v}
\end{bmatrix}
-
\begin{bmatrix}
O & I \\
B & O
\end{bmatrix}
\begin{bmatrix}
\mathbf{u} \\
\mathbf{v}
\end{bmatrix}
-
\begin{bmatrix}
0 \\
\beta \mathbf{u}^2
\end{bmatrix} = \mathbf{0},\quad
B = \frac{1}{\Delta x^2}
\begin{bmatrix}
1 & 0 &        &        & 0 \\
1  & -2 & 1      &        &   \\
   & \ddots & \ddots & \ddots &   \\
   &        & 1      & -2 & 1 \\
0  &        &        & 0 & 1
\end{bmatrix},
\]
where the operation $\mathbf{u}^2$ should be understood componentwise; see Fig.~\ref{fig:solution-of-wave} for a depiction of the solution of~\eqref{eq:nonlinearwaeve} with $\Delta x = \nicefrac{1}{127}$ and $h = \nicefrac{1}{1270}$ obtained using the symmetric Gauss scheme with $s = 4$ stages.
To have both a larger problem size and moving towards the direction of a more realistic problem, we consider the same test illustrated in Fig.~\ref{fig:solution-of-wave} discretized on a finer spatial mesh with $\Delta x = \nicefrac{1}{1023}$ and set the time-step size to $h = \nicefrac{\Delta x}{10}$. The simulation is carried out over $100$ time-steps to collect performance data, ensuring that all simulations reach the same final time independently of $s$. We employ the Gauss integration scheme with varying numbers of stages $s$, and initially focus on two metrics: the number of simplified Newton iterations required per time-step, and the average dimension of the Krylov subspace used in solving the correction Sylvester matrix equations. 
\begin{figure}[htbp]
    \centering
    \begin{subfigure}{0.3\columnwidth}\centering
    \input{iter_s5}
    \caption{5 stages}
    \end{subfigure}\hspace{1.2em}
    \begin{subfigure}{0.3\columnwidth}\centering
    \input{iter_s15}
    \caption{15 stages}
    \end{subfigure}\hspace{-0.5em}
    \begin{subfigure}{0.3\columnwidth}\centering
    \input{iter_s25}
    \caption{25 stages}
    \end{subfigure}
    \caption{Iteration count for the simplified Newton method and average dimension of the extended rational Krylov subspace for the correction matrix equation~\eqref{eq:sylvester_for_newton}.}
    \label{fig:iter_nonlin}
\end{figure}
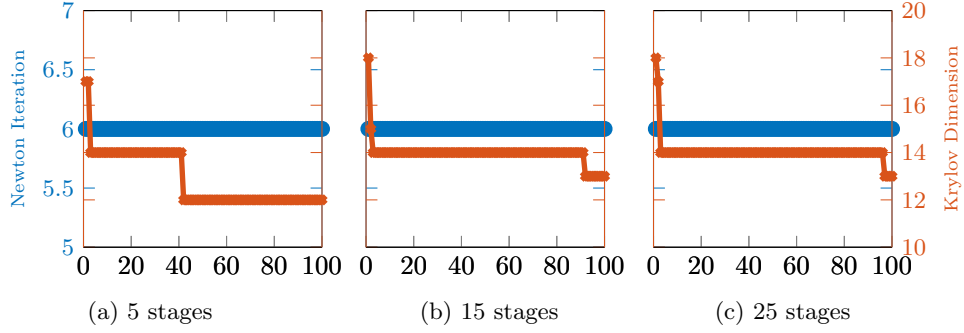
In particular, we utilize a variant of Algorithm~\ref{alg:oneside_sylvester} that employs the extended Krylov subspace method~\cite{MR2318706}, where the subspace is expanded by two basis vectors at each iteration. As shown in Fig.~\ref{fig:iter_nonlin}, the number of simplified Newton iterations—targeting a convergence tolerance of $10^{-10}$—remains stable across time-steps and is largely unaffected by the number of Gauss stages. The average dimension of the Krylov subspace exhibits some variability, but the fluctuations remain well-contained throughout the simulation.
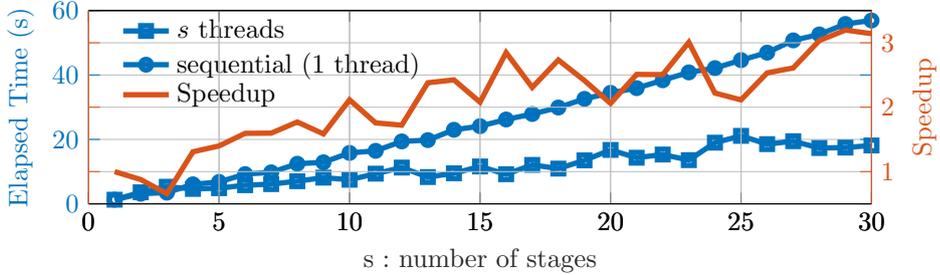
\begin{figure}[htbp]
    \centering
    \input{nonlinearperformance}
    
    \caption{Performance in terms of elapsed times and speedup of the threaded version against the sequential implementation for varying number of threads.}
    \label{fig:nonlinear_performances}
\end{figure}
Fig.~\ref{fig:nonlinear_performances} illustrates the performance comparison between the threaded and sequential implementations of the method as a function of the number of stages which is equal to the number of threads. The left axis reports the total elapsed time measured in seconds, while the right axis shows the speedup.  As expected, the execution time increases with the number of stages due to the higher computational cost per time-step. However, the threaded implementation consistently outperforms the sequential counterpart, achieving substantial speedup, particularly in the mid-to-high stage range. Notably, a speedup close to $3\times$ is observed for several stage counts, indicating effective parallel scalability of the threaded solver. The variability in speedup at higher stage counts may be attributed to increasing synchronization overhead or memory contention effects in the parallel execution. For low stage counts, however, the threaded implementation shows no significant performance improvement. This is primarily due to the limited amount of parallel work available in these configurations, where the overhead of thread management and synchronization outweighs the computational benefits. From a weak-scaling perspective, the behavior of the threaded implementation as \( s \) increases is encouraging:
to an increase from $2N$ to $30N$ corresponds an increase in time of $5\times$. Although the total work per time-step grows with \( s \), the threaded runtime exhibits only a moderate increase in elapsed time, suggesting that the solver sustains parallel efficiency even as the computational load per step increases. %

\section{Conclusions and future extensions}\label{sec:conclusions}

We have shown an approach for building stage-parallel IRK methods for high stage counts together with a feasible implementation in a multithreaded environment. The numerical results on two classical benchmark problems are promising and highlight the applicability of the method to more realistic settings. Among the future plans, there is to investigate in detail the construction of iterative methods for the solution of block diagonal systems that arise from the procedure of fast diagonalization of the equations for the stages. %
Moreover, we aim to extend the construction in two directions. First to encompass extension of Runge--Kutta methods for conservative problems which can be expressed in an analogous matrix format~\cite{BrugnanoBook}, secondly towards cases where, beyond handling the stages of the Runge–Kutta method, we also solve all time-steps simultaneously.
On the implementation side, we plan to incorporate this strategy into the \textit{Parallel Sparse Computation Toolkit}~\cite{DAmbra2023} to fully leverage a parallel computing environment. This approach enables the distribution of matrices across multiple processes and allows for the overlapping of stage computations, thereby enhancing overall performance. %

\backmatter

\bmhead{Acknowledgements}
FD acknowledges the MUR Excellence Department Project awarded to the Department of Mathematics, University of Pisa, CUP I57G22000700001.  MM acknowledges the MUR Excellence Department Project MatMod@TOV awarded to the Department of Mathematics, University of Rome Tor Vergata, CUP E83C23000330006. The research of FD was  partially granted by the Italian Ministry of University and Research (MUR) through the PRIN 2022 ``MOLE: Manifold constrained Optimization and LEarning'',  code: 2022ZK5ME7 MUR D.D. financing decree n. 20428 of November 6th, 2024 (CUP B53C24006410006). This project has received funding from the European High Performance Computing Joint Undertaking under grant agreement No. 101172493. Both authors are members of INdAM-GNCS and have been partially financed by the INdAM-GNCS Project CUP E53C24001950001. 

The authors thank the reviewers for their suggestion in improving the presentation of the paper.

\section*{Declarations}

\bmhead{Conflict of interest} The authors declare that they have no conflict of interest.
\bmhead{Data availability} Data sharing is not applicable to this article as no new datasets were generated
during the current study.
\bmhead{Code availability} The code discussed here is implemented in the Julia package available at the GitHub repository \href{https://github.com/Cirdans-Home/SP\_IRK.jl}{Cirdans-Home/SP\_IRK.jl}.
\bmhead{Author contribution} All authors contributed equally to the paper.

\bibliography{bibliography}

\end{document}

%% file: evcondition.tex
\begin{tikzpicture}

\begin{axis}[%
width=0.8\columnwidth,
height=1.2in,
at={(0.758in,0.423in)},
scale only axis,
xmin=2,
xmax=30,
xtick={ 2,  3,  4,  5,  6,  7,  8,  9, 10, 11, 12, 13, 14, 15, 16, 17, 18, 19, 20, 21, 22, 23, 24, 25, 26, 27, 28, 29, 30},
ytick={1,1e8,1e16},
xticklabel style={rotate=45},
xlabel style={font=\color{white!15!black}},
xlabel={s},
ymode=log,
ymin=1,
ymax=2.85754566935025e+16,
yminorticks=true,
axis background/.style={fill=white},
xmajorgrids,
ymajorgrids,
yminorgrids,
legend style={at={(0.03,0.97)}, anchor=north west, legend cell align=left, align=left, draw=white!15!black}
]
\addplot [color=red, line width=2.0pt, mark=square, mark options={solid, red}]
  table[row sep=crcr]{%
1	0\\
2	1\\
3	1.27855015329644\\
4	1.42726787499268\\
5	1.60562927331117\\
6	1.76798816128677\\
7	1.93211098829681\\
8	2.07841311146843\\
9	2.22810700381127\\
10	2.36359887169907\\
11	2.49798320487261\\
12	2.62347892667392\\
13	2.74696320475702\\
14	2.86305225182051\\
15	2.97777882608408\\
16	3.0861204516393\\
17	3.19315496340658\\
18	3.29503205530519\\
19	3.39543011979212\\
20	3.49163934320629\\
21	3.58648299035316\\
22	3.67765511674192\\
23	3.76774349562767\\
24	3.85452942013125\\
25	3.94042271039804\\
26	4.02335142167481\\
27	4.10550258961752\\
28	4.18500218683517\\
29	4.26377659971082\\
30	4.34018848975515\\
};
\addlegendentry{$\kappa_2(V)$}

\addplot [color=blue, line width=2.0pt, mark=+, mark options={solid, blue}]
  table[row sep=crcr]{%
1	0\\
2	3.73205080756888\\
3	12.8375509342729\\
4	45.5849425032636\\
5	165.004190977923\\
6	602.653866227601\\
7	2211.85081004009\\
8	8143.57238643118\\
9	30052.1252797002\\
10	111100.692003655\\
11	411331.64458242\\
12	1524732.48390779\\
13	5657683.70570603\\
14	21011750.4746833\\
15	78093241.1916797\\
16	290437443.977983\\
17	1080808122.08134\\
18	4024171883.29884\\
19	14990486472.9999\\
20	55866135332.3218\\
21	208285967023.266\\
22	776844787854.279\\
23	2898402118862.3\\
24	10817693852016.3\\
25	40382079473186.1\\
26	150782911992083\\
27	566155538393591\\
28	2.12458022781876e+15\\
29	8.10894173935367e+15\\
30	2.85754566935025e+16\\
};
\addlegendentry{$\kappa_2(X)$}

\end{axis}

\end{tikzpicture}%

%% file: iteration_time.tex
\definecolor{mycolor1}{rgb}{0.00000,0.44706,0.74118}%
\definecolor{mycolor2}{rgb}{0.85098,0.32549,0.09804}%
\begin{tikzpicture}

\begin{axis}[%
width=0.75\columnwidth,
height=0.174\columnwidth,
at={(0\columnwidth,0\columnwidth)},
scale only axis,
unbounded coords=jump,
xmin=0.5,
xmax=30.5,
xtick={1,2,3,4,5,6,7,8,9,10,11,12,13,14,15,16,17,18,19,20,21,22,23,24,25,26,27,28,29,30},
xticklabel style = {rotate=90,font=\small},
xlabel style={font=\color{white!15!black}},
xlabel={s : number of stages / threads},
separate axis lines,
every outer y axis line/.append style={mycolor1},
every y tick label/.append style={font=\color{mycolor1}},
every y tick/.append style={mycolor1},
ymin=6,
ymax=80,
ylabel style={font=\color{mycolor1}},
ylabel={Elapsed Time (s)},
]
\addplot [color=mycolor1, line width=2.0pt, forget plot]
  table[row sep=crcr]{%
1	6.114\\
2	18.86\\
3	14.199\\
4	24.784\\
5	19.841\\
6	26.067\\
7	23.197\\
8	29.239\\
9	27.414\\
10	31.814\\
11	28.974\\
12	35.414\\
13	33.674\\
14	37.169\\
15	36.698\\
16	39.83\\
17	58.749\\
18	58.152\\
19	72.197\\
20	63.276\\
21	67.598\\
22	69.087\\
23	63.172\\
24	61.682\\
25	72.732\\
26	67.287\\
27	65.405\\
28	71.72\\
29	62.911\\
30	70.11\\
};
\end{axis}
\begin{axis}[%
width=0.75\columnwidth,
height=0.174\columnwidth,
at={(0\columnwidth,0\columnwidth)},
scale only axis,
unbounded coords=jump,
xmin=0.5,
xmax=30.5,
xtick={},
xticklabels={},
separate axis lines,
every outer y axis line/.append style={mycolor2},
every y tick label/.append style={font=\color{mycolor2}},
every y tick/.append style={mycolor2},
ymin=6,
ymax=101,
ylabel style={font=\color{mycolor2}},
ylabel={Krylov Iteration},
axis background/.style={fill=none},
yticklabel pos=right
]
\addplot [color=black, dashed, forget plot]
  table[row sep=crcr]{%
1	16\\
1	17\\
};
\addplot [color=black, dashed, forget plot]
  table[row sep=crcr]{%
2	26\\
2	27\\
};
\addplot [color=black, dashed, forget plot]
  table[row sep=crcr]{%
3	25\\
3	26\\
};
\addplot [color=black, dashed, forget plot]
  table[row sep=crcr]{%
4	33\\
4	36\\
};
\addplot [color=black, dashed, forget plot]
  table[row sep=crcr]{%
5	34\\
5	37\\
};
\addplot [color=black, dashed, forget plot]
  table[row sep=crcr]{%
6	39\\
6	43\\
};
\addplot [color=black, dashed, forget plot]
  table[row sep=crcr]{%
7	42\\
7	47\\
};
\addplot [color=black, dashed, forget plot]
  table[row sep=crcr]{%
8	44\\
8	48\\
};
\addplot [color=black, dashed, forget plot]
  table[row sep=crcr]{%
9	47\\
9	49\\
};
\addplot [color=black, dashed, forget plot]
  table[row sep=crcr]{%
10	50\\
10	51\\
};
\addplot [color=black, dashed, forget plot]
  table[row sep=crcr]{%
11	53\\
11	54\\
};
\addplot [color=black, dashed, forget plot]
  table[row sep=crcr]{%
12	55.5\\
12	59\\
};
\addplot [color=black, dashed, forget plot]
  table[row sep=crcr]{%
13	58.5\\
13	62\\
};
\addplot [color=black, dashed, forget plot]
  table[row sep=crcr]{%
14	61\\
14	66\\
};
\addplot [color=black, dashed, forget plot]
  table[row sep=crcr]{%
15	64\\
15	69\\
};
\addplot [color=black, dashed, forget plot]
  table[row sep=crcr]{%
16	66\\
16	72\\
};
\addplot [color=black, dashed, forget plot]
  table[row sep=crcr]{%
17	69\\
17	75\\
};
\addplot [color=black, dashed, forget plot]
  table[row sep=crcr]{%
18	71\\
18	77\\
};
\addplot [color=black, dashed, forget plot]
  table[row sep=crcr]{%
19	73\\
19	81\\
};
\addplot [color=black, dashed, forget plot]
  table[row sep=crcr]{%
20	75\\
20	83\\
};
\addplot [color=black, dashed, forget plot]
  table[row sep=crcr]{%
21	77\\
21	86\\
};
\addplot [color=black, dashed, forget plot]
  table[row sep=crcr]{%
22	79\\
22	88\\
};
\addplot [color=black, dashed, forget plot]
  table[row sep=crcr]{%
23	81\\
23	91\\
};
\addplot [color=black, dashed, forget plot]
  table[row sep=crcr]{%
24	83\\
24	93\\
};
\addplot [color=black, dashed, forget plot]
  table[row sep=crcr]{%
25	85\\
25	96\\
};
\addplot [color=black, dashed, forget plot]
  table[row sep=crcr]{%
26	87\\
26	98\\
};
\addplot [color=black, dashed, forget plot]
  table[row sep=crcr]{%
27	88.5\\
27	100\\
};
\addplot [color=black, dashed, forget plot]
  table[row sep=crcr]{%
28	90\\
28	100\\
};
\addplot [color=black, dashed, forget plot]
  table[row sep=crcr]{%
29	92\\
29	100\\
};
\addplot [color=black, dashed, forget plot]
  table[row sep=crcr]{%
30	94\\
30	100\\
};
\addplot [color=black, dashed, forget plot]
  table[row sep=crcr]{%
1	14\\
1	15\\
};
\addplot [color=black, dashed, forget plot]
  table[row sep=crcr]{%
2	22\\
2	24\\
};
\addplot [color=black, dashed, forget plot]
  table[row sep=crcr]{%
3	23\\
3	24\\
};
\addplot [color=black, dashed, forget plot]
  table[row sep=crcr]{%
4	28\\
4	30\\
};
\addplot [color=black, dashed, forget plot]
  table[row sep=crcr]{%
5	29\\
5	32\\
};
\addplot [color=black, dashed, forget plot]
  table[row sep=crcr]{%
6	33\\
6	36\\
};
\addplot [color=black, dashed, forget plot]
  table[row sep=crcr]{%
7	35\\
7	38\\
};
\addplot [color=black, dashed, forget plot]
  table[row sep=crcr]{%
8	39\\
8	41\\
};
\addplot [color=black, dashed, forget plot]
  table[row sep=crcr]{%
9	42\\
9	44\\
};
\addplot [color=black, dashed, forget plot]
  table[row sep=crcr]{%
10	44\\
10	47\\
};
\addplot [color=black, dashed, forget plot]
  table[row sep=crcr]{%
11	47\\
11	50\\
};
\addplot [color=black, dashed, forget plot]
  table[row sep=crcr]{%
12	49\\
12	52\\
};
\addplot [color=black, dashed, forget plot]
  table[row sep=crcr]{%
13	51\\
13	55\\
};
\addplot [color=black, dashed, forget plot]
  table[row sep=crcr]{%
14	52\\
14	57\\
};
\addplot [color=black, dashed, forget plot]
  table[row sep=crcr]{%
15	55\\
15	59.5\\
};
\addplot [color=black, dashed, forget plot]
  table[row sep=crcr]{%
16	55\\
16	61\\
};
\addplot [color=black, dashed, forget plot]
  table[row sep=crcr]{%
17	58\\
17	64\\
};
\addplot [color=black, dashed, forget plot]
  table[row sep=crcr]{%
18	58\\
18	65\\
};
\addplot [color=black, dashed, forget plot]
  table[row sep=crcr]{%
19	61\\
19	67\\
};
\addplot [color=black, dashed, forget plot]
  table[row sep=crcr]{%
20	61\\
20	69\\
};
\addplot [color=black, dashed, forget plot]
  table[row sep=crcr]{%
21	64\\
21	71\\
};
\addplot [color=black, dashed, forget plot]
  table[row sep=crcr]{%
22	64\\
22	72.5\\
};
\addplot [color=black, dashed, forget plot]
  table[row sep=crcr]{%
23	66\\
23	74\\
};
\addplot [color=black, dashed, forget plot]
  table[row sep=crcr]{%
24	67\\
24	75.5\\
};
\addplot [color=black, dashed, forget plot]
  table[row sep=crcr]{%
25	68\\
25	76.5\\
};
\addplot [color=black, dashed, forget plot]
  table[row sep=crcr]{%
26	69\\
26	77.5\\
};
\addplot [color=black, dashed, forget plot]
  table[row sep=crcr]{%
27	70\\
27	79.5\\
};
\addplot [color=black, dashed, forget plot]
  table[row sep=crcr]{%
28	71\\
28	80\\
};
\addplot [color=black, dashed, forget plot]
  table[row sep=crcr]{%
29	72\\
29	82\\
};
\addplot [color=black, dashed, forget plot]
  table[row sep=crcr]{%
30	72\\
30	82.5\\
};
\addplot [color=black, forget plot]
  table[row sep=crcr]{%
0.875	17\\
1.125	17\\
};
\addplot [color=black, forget plot]
  table[row sep=crcr]{%
1.875	27\\
2.125	27\\
};
\addplot [color=black, forget plot]
  table[row sep=crcr]{%
2.875	26\\
3.125	26\\
};
\addplot [color=black, forget plot]
  table[row sep=crcr]{%
3.875	36\\
4.125	36\\
};
\addplot [color=black, forget plot]
  table[row sep=crcr]{%
4.875	37\\
5.125	37\\
};
\addplot [color=black, forget plot]
  table[row sep=crcr]{%
5.875	43\\
6.125	43\\
};
\addplot [color=black, forget plot]
  table[row sep=crcr]{%
6.875	47\\
7.125	47\\
};
\addplot [color=black, forget plot]
  table[row sep=crcr]{%
7.875	48\\
8.125	48\\
};
\addplot [color=black, forget plot]
  table[row sep=crcr]{%
8.875	49\\
9.125	49\\
};
\addplot [color=black, forget plot]
  table[row sep=crcr]{%
9.875	51\\
10.125	51\\
};
\addplot [color=black, forget plot]
  table[row sep=crcr]{%
10.875	54\\
11.125	54\\
};
\addplot [color=black, forget plot]
  table[row sep=crcr]{%
11.875	59\\
12.125	59\\
};
\addplot [color=black, forget plot]
  table[row sep=crcr]{%
12.875	62\\
13.125	62\\
};
\addplot [color=black, forget plot]
  table[row sep=crcr]{%
13.875	66\\
14.125	66\\
};
\addplot [color=black, forget plot]
  table[row sep=crcr]{%
14.875	69\\
15.125	69\\
};
\addplot [color=black, forget plot]
  table[row sep=crcr]{%
15.875	72\\
16.125	72\\
};
\addplot [color=black, forget plot]
  table[row sep=crcr]{%
16.875	75\\
17.125	75\\
};
\addplot [color=black, forget plot]
  table[row sep=crcr]{%
17.875	77\\
18.125	77\\
};
\addplot [color=black, forget plot]
  table[row sep=crcr]{%
18.875	81\\
19.125	81\\
};
\addplot [color=black, forget plot]
  table[row sep=crcr]{%
19.875	83\\
20.125	83\\
};
\addplot [color=black, forget plot]
  table[row sep=crcr]{%
20.875	86\\
21.125	86\\
};
\addplot [color=black, forget plot]
  table[row sep=crcr]{%
21.875	88\\
22.125	88\\
};
\addplot [color=black, forget plot]
  table[row sep=crcr]{%
22.875	91\\
23.125	91\\
};
\addplot [color=black, forget plot]
  table[row sep=crcr]{%
23.875	93\\
24.125	93\\
};
\addplot [color=black, forget plot]
  table[row sep=crcr]{%
24.875	96\\
25.125	96\\
};
\addplot [color=black, forget plot]
  table[row sep=crcr]{%
25.875	98\\
26.125	98\\
};
\addplot [color=black, forget plot]
  table[row sep=crcr]{%
26.875	100\\
27.125	100\\
};
\addplot [color=black, forget plot]
  table[row sep=crcr]{%
27.875	100\\
28.125	100\\
};
\addplot [color=black, forget plot]
  table[row sep=crcr]{%
28.875	100\\
29.125	100\\
};
\addplot [color=black, forget plot]
  table[row sep=crcr]{%
29.875	100\\
30.125	100\\
};
\addplot [color=black, forget plot]
  table[row sep=crcr]{%
0.875	14\\
1.125	14\\
};
\addplot [color=black, forget plot]
  table[row sep=crcr]{%
1.875	22\\
2.125	22\\
};
\addplot [color=black, forget plot]
  table[row sep=crcr]{%
2.875	23\\
3.125	23\\
};
\addplot [color=black, forget plot]
  table[row sep=crcr]{%
3.875	28\\
4.125	28\\
};
\addplot [color=black, forget plot]
  table[row sep=crcr]{%
4.875	29\\
5.125	29\\
};
\addplot [color=black, forget plot]
  table[row sep=crcr]{%
5.875	33\\
6.125	33\\
};
\addplot [color=black, forget plot]
  table[row sep=crcr]{%
6.875	35\\
7.125	35\\
};
\addplot [color=black, forget plot]
  table[row sep=crcr]{%
7.875	39\\
8.125	39\\
};
\addplot [color=black, forget plot]
  table[row sep=crcr]{%
8.875	42\\
9.125	42\\
};
\addplot [color=black, forget plot]
  table[row sep=crcr]{%
9.875	44\\
10.125	44\\
};
\addplot [color=black, forget plot]
  table[row sep=crcr]{%
10.875	47\\
11.125	47\\
};
\addplot [color=black, forget plot]
  table[row sep=crcr]{%
11.875	49\\
12.125	49\\
};
\addplot [color=black, forget plot]
  table[row sep=crcr]{%
12.875	51\\
13.125	51\\
};
\addplot [color=black, forget plot]
  table[row sep=crcr]{%
13.875	52\\
14.125	52\\
};
\addplot [color=black, forget plot]
  table[row sep=crcr]{%
14.875	55\\
15.125	55\\
};
\addplot [color=black, forget plot]
  table[row sep=crcr]{%
15.875	55\\
16.125	55\\
};
\addplot [color=black, forget plot]
  table[row sep=crcr]{%
16.875	58\\
17.125	58\\
};
\addplot [color=black, forget plot]
  table[row sep=crcr]{%
17.875	58\\
18.125	58\\
};
\addplot [color=black, forget plot]
  table[row sep=crcr]{%
18.875	61\\
19.125	61\\
};
\addplot [color=black, forget plot]
  table[row sep=crcr]{%
19.875	61\\
20.125	61\\
};
\addplot [color=black, forget plot]
  table[row sep=crcr]{%
20.875	64\\
21.125	64\\
};
\addplot [color=black, forget plot]
  table[row sep=crcr]{%
21.875	64\\
22.125	64\\
};
\addplot [color=black, forget plot]
  table[row sep=crcr]{%
22.875	66\\
23.125	66\\
};
\addplot [color=black, forget plot]
  table[row sep=crcr]{%
23.875	67\\
24.125	67\\
};
\addplot [color=black, forget plot]
  table[row sep=crcr]{%
24.875	68\\
25.125	68\\
};
\addplot [color=black, forget plot]
  table[row sep=crcr]{%
25.875	69\\
26.125	69\\
};
\addplot [color=black, forget plot]
  table[row sep=crcr]{%
26.875	70\\
27.125	70\\
};
\addplot [color=black, forget plot]
  table[row sep=crcr]{%
27.875	71\\
28.125	71\\
};
\addplot [color=black, forget plot]
  table[row sep=crcr]{%
28.875	72\\
29.125	72\\
};
\addplot [color=black, forget plot]
  table[row sep=crcr]{%
29.875	72\\
30.125	72\\
};
\addplot [color=blue, forget plot]
  table[row sep=crcr]{%
0.75	15\\
0.75	16\\
1.25	16\\
1.25	15\\
0.75	15\\
};
\addplot [color=blue, forget plot]
  table[row sep=crcr]{%
1.75	24\\
1.75	26\\
2.25	26\\
2.25	24\\
1.75	24\\
};
\addplot [color=blue, forget plot]
  table[row sep=crcr]{%
2.75	24\\
2.75	25\\
3.25	25\\
3.25	24\\
2.75	24\\
};
\addplot [color=blue, forget plot]
  table[row sep=crcr]{%
3.75	30\\
3.75	33\\
4.25	33\\
4.25	30\\
3.75	30\\
};
\addplot [color=blue, forget plot]
  table[row sep=crcr]{%
4.75	32\\
4.75	34\\
5.25	34\\
5.25	32\\
4.75	32\\
};
\addplot [color=blue, forget plot]
  table[row sep=crcr]{%
5.75	36\\
5.75	39\\
6.25	39\\
6.25	36\\
5.75	36\\
};
\addplot [color=blue, forget plot]
  table[row sep=crcr]{%
6.75	38\\
6.75	42\\
7.25	42\\
7.25	38\\
6.75	38\\
};
\addplot [color=blue, forget plot]
  table[row sep=crcr]{%
7.75	41\\
7.75	44\\
8.25	44\\
8.25	41\\
7.75	41\\
};
\addplot [color=blue, forget plot]
  table[row sep=crcr]{%
8.75	44\\
8.75	47\\
9.25	47\\
9.25	44\\
8.75	44\\
};
\addplot [color=blue, forget plot]
  table[row sep=crcr]{%
9.75	47\\
9.75	50\\
10.25	50\\
10.25	47\\
9.75	47\\
};
\addplot [color=blue, forget plot]
  table[row sep=crcr]{%
10.75	50\\
10.75	53\\
11.25	53\\
11.25	50\\
10.75	50\\
};
\addplot [color=blue, forget plot]
  table[row sep=crcr]{%
11.75	52\\
11.75	55.5\\
12.25	55.5\\
12.25	52\\
11.75	52\\
};
\addplot [color=blue, forget plot]
  table[row sep=crcr]{%
12.75	55\\
12.75	58.5\\
13.25	58.5\\
13.25	55\\
12.75	55\\
};
\addplot [color=blue, forget plot]
  table[row sep=crcr]{%
13.75	57\\
13.75	61\\
14.25	61\\
14.25	57\\
13.75	57\\
};
\addplot [color=blue, forget plot]
  table[row sep=crcr]{%
14.75	59.5\\
14.75	64\\
15.25	64\\
15.25	59.5\\
14.75	59.5\\
};
\addplot [color=blue, forget plot]
  table[row sep=crcr]{%
15.75	61\\
15.75	66\\
16.25	66\\
16.25	61\\
15.75	61\\
};
\addplot [color=blue, forget plot]
  table[row sep=crcr]{%
16.75	64\\
16.75	69\\
17.25	69\\
17.25	64\\
16.75	64\\
};
\addplot [color=blue, forget plot]
  table[row sep=crcr]{%
17.75	65\\
17.75	71\\
18.25	71\\
18.25	65\\
17.75	65\\
};
\addplot [color=blue, forget plot]
  table[row sep=crcr]{%
18.75	67\\
18.75	73\\
19.25	73\\
19.25	67\\
18.75	67\\
};
\addplot [color=blue, forget plot]
  table[row sep=crcr]{%
19.75	69\\
19.75	75\\
20.25	75\\
20.25	69\\
19.75	69\\
};
\addplot [color=blue, forget plot]
  table[row sep=crcr]{%
20.75	71\\
20.75	77\\
21.25	77\\
21.25	71\\
20.75	71\\
};
\addplot [color=blue, forget plot]
  table[row sep=crcr]{%
21.75	72.5\\
21.75	79\\
22.25	79\\
22.25	72.5\\
21.75	72.5\\
};
\addplot [color=blue, forget plot]
  table[row sep=crcr]{%
22.75	74\\
22.75	81\\
23.25	81\\
23.25	74\\
22.75	74\\
};
\addplot [color=blue, forget plot]
  table[row sep=crcr]{%
23.75	75.5\\
23.75	83\\
24.25	83\\
24.25	75.5\\
23.75	75.5\\
};
\addplot [color=blue, forget plot]
  table[row sep=crcr]{%
24.75	76.5\\
24.75	85\\
25.25	85\\
25.25	76.5\\
24.75	76.5\\
};
\addplot [color=blue, forget plot]
  table[row sep=crcr]{%
25.75	77.5\\
25.75	87\\
26.25	87\\
26.25	77.5\\
25.75	77.5\\
};
\addplot [color=blue, forget plot]
  table[row sep=crcr]{%
26.75	79.5\\
26.75	88.5\\
27.25	88.5\\
27.25	79.5\\
26.75	79.5\\
};
\addplot [color=blue, forget plot]
  table[row sep=crcr]{%
27.75	80\\
27.75	90\\
28.25	90\\
28.25	80\\
27.75	80\\
};
\addplot [color=blue, forget plot]
  table[row sep=crcr]{%
28.75	82\\
28.75	92\\
29.25	92\\
29.25	82\\
28.75	82\\
};
\addplot [color=blue, forget plot]
  table[row sep=crcr]{%
29.75	82.5\\
29.75	94\\
30.25	94\\
30.25	82.5\\
29.75	82.5\\
};
\addplot [color=red, forget plot]
  table[row sep=crcr]{%
0.75	15\\
1.25	15\\
};
\addplot [color=red, forget plot]
  table[row sep=crcr]{%
1.75	25\\
2.25	25\\
};
\addplot [color=red, forget plot]
  table[row sep=crcr]{%
2.75	25\\
3.25	25\\
};
\addplot [color=red, forget plot]
  table[row sep=crcr]{%
3.75	32\\
4.25	32\\
};
\addplot [color=red, forget plot]
  table[row sep=crcr]{%
4.75	34\\
5.25	34\\
};
\addplot [color=red, forget plot]
  table[row sep=crcr]{%
5.75	39\\
6.25	39\\
};
\addplot [color=red, forget plot]
  table[row sep=crcr]{%
6.75	42\\
7.25	42\\
};
\addplot [color=red, forget plot]
  table[row sep=crcr]{%
7.75	43\\
8.25	43\\
};
\addplot [color=red, forget plot]
  table[row sep=crcr]{%
8.75	46\\
9.25	46\\
};
\addplot [color=red, forget plot]
  table[row sep=crcr]{%
9.75	49\\
10.25	49\\
};
\addplot [color=red, forget plot]
  table[row sep=crcr]{%
10.75	52\\
11.25	52\\
};
\addplot [color=red, forget plot]
  table[row sep=crcr]{%
11.75	55\\
12.25	55\\
};
\addplot [color=red, forget plot]
  table[row sep=crcr]{%
12.75	58\\
13.25	58\\
};
\addplot [color=red, forget plot]
  table[row sep=crcr]{%
13.75	60\\
14.25	60\\
};
\addplot [color=red, forget plot]
  table[row sep=crcr]{%
14.75	63\\
15.25	63\\
};
\addplot [color=red, forget plot]
  table[row sep=crcr]{%
15.75	65\\
16.25	65\\
};
\addplot [color=red, forget plot]
  table[row sep=crcr]{%
16.75	68\\
17.25	68\\
};
\addplot [color=red, forget plot]
  table[row sep=crcr]{%
17.75	70\\
18.25	70\\
};
\addplot [color=red, forget plot]
  table[row sep=crcr]{%
18.75	72\\
19.25	72\\
};
\addplot [color=red, forget plot]
  table[row sep=crcr]{%
19.75	74\\
20.25	74\\
};
\addplot [color=red, forget plot]
  table[row sep=crcr]{%
20.75	76\\
21.25	76\\
};
\addplot [color=red, forget plot]
  table[row sep=crcr]{%
21.75	78\\
22.25	78\\
};
\addplot [color=red, forget plot]
  table[row sep=crcr]{%
22.75	80\\
23.25	80\\
};
\addplot [color=red, forget plot]
  table[row sep=crcr]{%
23.75	81.5\\
24.25	81.5\\
};
\addplot [color=red, forget plot]
  table[row sep=crcr]{%
24.75	83.5\\
25.25	83.5\\
};
\addplot [color=red, forget plot]
  table[row sep=crcr]{%
25.75	85\\
26.25	85\\
};
\addplot [color=red, forget plot]
  table[row sep=crcr]{%
26.75	87\\
27.25	87\\
};
\addplot [color=red, forget plot]
  table[row sep=crcr]{%
27.75	88\\
28.25	88\\
};
\addplot [color=red, forget plot]
  table[row sep=crcr]{%
28.75	90\\
29.25	90\\
};
\addplot [color=red, forget plot]
  table[row sep=crcr]{%
29.75	91.5\\
30.25	91.5\\
};
\addplot [color=black, only marks, mark=+, mark options={solid, draw=red}, forget plot]
  table[row sep=crcr]{%
1	13\\
1	18\\
1	18\\
1	18\\
1	19\\
1	19\\
1	20\\
1	21\\
};
\addplot [color=black, only marks, mark=+, mark options={solid, draw=red}, forget plot]
  table[row sep=crcr]{%
nan	nan\\
};
\addplot [color=black, only marks, mark=+, mark options={solid, draw=red}, forget plot]
  table[row sep=crcr]{%
3	21\\
3	22\\
3	27\\
3	28\\
3	30\\
3	32\\
};
\end{axis}
\end{tikzpicture}%

%% file: iteration_time_radau.tex
\definecolor{mycolor1}{rgb}{0.00000,0.44706,0.74118}%
\definecolor{mycolor2}{rgb}{0.85098,0.32549,0.09804}%
\begin{tikzpicture}

\begin{axis}[%
width=0.75\columnwidth,
height=0.174\columnwidth,
at={(0\columnwidth,0\columnwidth)},
scale only axis,
unbounded coords=jump,
xmin=2,
xmax=30,
xtick={1,2,3,4,5,6,7,8,9,10,11,12,13,14,15,16,17,18,19,20,21,22,23,24,25,26,27,28,29,30},
xticklabel style = {rotate=90,font=\small},
xlabel style={font=\color{white!15!black}},
xlabel={s : number of stages / threads},
separate axis lines,
every outer y axis line/.append style={mycolor1},
every y tick label/.append style={font=\color{mycolor1}},
every y tick/.append style={mycolor2},
ymin=2.75,
ymax=70.00,
ylabel style={font=\color{mycolor1}},
ylabel={Elapsed Time (s)},
]
\addplot [color=mycolor1, line width=2.0pt, forget plot]
  table[row sep=crcr]{%
2	16.488\\
3	13.446\\
4	19.768\\
5	17.366\\
6	23.995\\
7	19.845\\
8	25.728\\
9	22.198\\
10	26.017\\
11	25.076\\
12	29.476\\
13	27.293\\
14	30.688\\
15	29.188\\
16	32.675\\
17	57.317\\
18	46.875\\
19	55.15\\
20	53.39\\
21	62.438\\
22	58.94\\
23	60.504\\
24	55.777\\
25	64.969\\
26	51.586\\
27	63.959\\
28	54.026\\
29	55.35\\
30	57.771\\
};
\end{axis}
\begin{axis}[%
width=0.75\columnwidth,
height=0.174\columnwidth,
at={(0\columnwidth,0\columnwidth)},
scale only axis,
unbounded coords=jump,
xmin=2,
xmax=30,
xtick={},
xticklabels={},
separate axis lines,
every outer y axis line/.append style={mycolor2},
every y tick label/.append style={font=\color{mycolor2}},
every y tick/.append style={mycolor2},
ymin=6,
ymax=70,
ylabel style={font=\color{mycolor2}},
ylabel={Krylov Iteration},
axis background/.style={fill=none},
yticklabel pos=right
]
\addplot [color=black, dashed, forget plot]
  table[row sep=crcr]{%
1	nan\\
1	nan\\
};
\addplot [color=black, dashed, forget plot]
  table[row sep=crcr]{%
2	8\\
2	11\\
};
\addplot [color=black, dashed, forget plot]
  table[row sep=crcr]{%
3	8\\
3	11\\
};
\addplot [color=black, dashed, forget plot]
  table[row sep=crcr]{%
4	8.5\\
4	12\\
};
\addplot [color=black, dashed, forget plot]
  table[row sep=crcr]{%
5	9\\
5	12\\
};
\addplot [color=black, dashed, forget plot]
  table[row sep=crcr]{%
6	10\\
6	11\\
};
\addplot [color=black, dashed, forget plot]
  table[row sep=crcr]{%
7	11\\
7	12\\
};
\addplot [color=black, dashed, forget plot]
  table[row sep=crcr]{%
8	13\\
8	14\\
};
\addplot [color=black, dashed, forget plot]
  table[row sep=crcr]{%
9	15\\
9	16\\
};
\addplot [color=black, dashed, forget plot]
  table[row sep=crcr]{%
10	17\\
10	18\\
};
\addplot [color=black, dashed, forget plot]
  table[row sep=crcr]{%
11	19\\
11	21\\
};
\addplot [color=black, dashed, forget plot]
  table[row sep=crcr]{%
12	20\\
12	22\\
};
\addplot [color=black, dashed, forget plot]
  table[row sep=crcr]{%
13	22\\
13	23\\
};
\addplot [color=black, dashed, forget plot]
  table[row sep=crcr]{%
14	23\\
14	26\\
};
\addplot [color=black, dashed, forget plot]
  table[row sep=crcr]{%
15	24\\
15	27\\
};
\addplot [color=black, dashed, forget plot]
  table[row sep=crcr]{%
16	25\\
16	27\\
};
\addplot [color=black, dashed, forget plot]
  table[row sep=crcr]{%
17	26\\
17	28\\
};
\addplot [color=black, dashed, forget plot]
  table[row sep=crcr]{%
18	28\\
18	32\\
};
\addplot [color=black, dashed, forget plot]
  table[row sep=crcr]{%
19	29\\
19	31\\
};
\addplot [color=black, dashed, forget plot]
  table[row sep=crcr]{%
20	29\\
20	31\\
};
\addplot [color=black, dashed, forget plot]
  table[row sep=crcr]{%
21	30\\
21	32\\
};
\addplot [color=black, dashed, forget plot]
  table[row sep=crcr]{%
22	32\\
22	38\\
};
\addplot [color=black, dashed, forget plot]
  table[row sep=crcr]{%
23	32\\
23	34\\
};
\addplot [color=black, dashed, forget plot]
  table[row sep=crcr]{%
24	33\\
24	38\\
};
\addplot [color=black, dashed, forget plot]
  table[row sep=crcr]{%
25	34.5\\
25	38\\
};
\addplot [color=black, dashed, forget plot]
  table[row sep=crcr]{%
26	36\\
26	42\\
};
\addplot [color=black, dashed, forget plot]
  table[row sep=crcr]{%
27	36\\
27	40\\
};
\addplot [color=black, dashed, forget plot]
  table[row sep=crcr]{%
28	37\\
28	43\\
};
\addplot [color=black, dashed, forget plot]
  table[row sep=crcr]{%
29	38\\
29	44\\
};
\addplot [color=black, dashed, forget plot]
  table[row sep=crcr]{%
30	39\\
30	43\\
};
\addplot [color=black, dashed, forget plot]
  table[row sep=crcr]{%
1	nan\\
1	nan\\
};
\addplot [color=black, dashed, forget plot]
  table[row sep=crcr]{%
2	5\\
2	6\\
};
\addplot [color=black, dashed, forget plot]
  table[row sep=crcr]{%
3	5\\
3	6\\
};
\addplot [color=black, dashed, forget plot]
  table[row sep=crcr]{%
4	6\\
4	6\\
};
\addplot [color=black, dashed, forget plot]
  table[row sep=crcr]{%
5	6\\
5	7\\
};
\addplot [color=black, dashed, forget plot]
  table[row sep=crcr]{%
6	8\\
6	9\\
};
\addplot [color=black, dashed, forget plot]
  table[row sep=crcr]{%
7	10\\
7	10\\
};
\addplot [color=black, dashed, forget plot]
  table[row sep=crcr]{%
8	11\\
8	12\\
};
\addplot [color=black, dashed, forget plot]
  table[row sep=crcr]{%
9	13\\
9	14\\
};
\addplot [color=black, dashed, forget plot]
  table[row sep=crcr]{%
10	15\\
10	16\\
};
\addplot [color=black, dashed, forget plot]
  table[row sep=crcr]{%
11	15\\
11	17\\
};
\addplot [color=black, dashed, forget plot]
  table[row sep=crcr]{%
12	16\\
12	18\\
};
\addplot [color=black, dashed, forget plot]
  table[row sep=crcr]{%
13	17\\
13	19.5\\
};
\addplot [color=black, dashed, forget plot]
  table[row sep=crcr]{%
14	18\\
14	20.5\\
};
\addplot [color=black, dashed, forget plot]
  table[row sep=crcr]{%
15	19\\
15	21.5\\
};
\addplot [color=black, dashed, forget plot]
  table[row sep=crcr]{%
16	20\\
16	23\\
};
\addplot [color=black, dashed, forget plot]
  table[row sep=crcr]{%
17	21\\
17	23.5\\
};
\addplot [color=black, dashed, forget plot]
  table[row sep=crcr]{%
18	22\\
18	24.5\\
};
\addplot [color=black, dashed, forget plot]
  table[row sep=crcr]{%
19	23\\
19	26\\
};
\addplot [color=black, dashed, forget plot]
  table[row sep=crcr]{%
20	23\\
20	26.5\\
};
\addplot [color=black, dashed, forget plot]
  table[row sep=crcr]{%
21	25\\
21	27\\
};
\addplot [color=black, dashed, forget plot]
  table[row sep=crcr]{%
22	25\\
22	28\\
};
\addplot [color=black, dashed, forget plot]
  table[row sep=crcr]{%
23	26\\
23	29\\
};
\addplot [color=black, dashed, forget plot]
  table[row sep=crcr]{%
24	26\\
24	29.5\\
};
\addplot [color=black, dashed, forget plot]
  table[row sep=crcr]{%
25	27\\
25	30.5\\
};
\addplot [color=black, dashed, forget plot]
  table[row sep=crcr]{%
26	28\\
26	31\\
};
\addplot [color=black, dashed, forget plot]
  table[row sep=crcr]{%
27	29\\
27	32\\
};
\addplot [color=black, dashed, forget plot]
  table[row sep=crcr]{%
28	29\\
28	33\\
};
\addplot [color=black, dashed, forget plot]
  table[row sep=crcr]{%
29	29\\
29	33.5\\
};
\addplot [color=black, dashed, forget plot]
  table[row sep=crcr]{%
30	31\\
30	34\\
};
\addplot [color=black, forget plot]
  table[row sep=crcr]{%
0.875	nan\\
1.125	nan\\
};
\addplot [color=black, forget plot]
  table[row sep=crcr]{%
1.875	11\\
2.125	11\\
};
\addplot [color=black, forget plot]
  table[row sep=crcr]{%
2.875	11\\
3.125	11\\
};
\addplot [color=black, forget plot]
  table[row sep=crcr]{%
3.875	12\\
4.125	12\\
};
\addplot [color=black, forget plot]
  table[row sep=crcr]{%
4.875	12\\
5.125	12\\
};
\addplot [color=black, forget plot]
  table[row sep=crcr]{%
5.875	11\\
6.125	11\\
};
\addplot [color=black, forget plot]
  table[row sep=crcr]{%
6.875	12\\
7.125	12\\
};
\addplot [color=black, forget plot]
  table[row sep=crcr]{%
7.875	14\\
8.125	14\\
};
\addplot [color=black, forget plot]
  table[row sep=crcr]{%
8.875	16\\
9.125	16\\
};
\addplot [color=black, forget plot]
  table[row sep=crcr]{%
9.875	18\\
10.125	18\\
};
\addplot [color=black, forget plot]
  table[row sep=crcr]{%
10.875	21\\
11.125	21\\
};
\addplot [color=black, forget plot]
  table[row sep=crcr]{%
11.875	22\\
12.125	22\\
};
\addplot [color=black, forget plot]
  table[row sep=crcr]{%
12.875	23\\
13.125	23\\
};
\addplot [color=black, forget plot]
  table[row sep=crcr]{%
13.875	26\\
14.125	26\\
};
\addplot [color=black, forget plot]
  table[row sep=crcr]{%
14.875	27\\
15.125	27\\
};
\addplot [color=black, forget plot]
  table[row sep=crcr]{%
15.875	27\\
16.125	27\\
};
\addplot [color=black, forget plot]
  table[row sep=crcr]{%
16.875	28\\
17.125	28\\
};
\addplot [color=black, forget plot]
  table[row sep=crcr]{%
17.875	32\\
18.125	32\\
};
\addplot [color=black, forget plot]
  table[row sep=crcr]{%
18.875	31\\
19.125	31\\
};
\addplot [color=black, forget plot]
  table[row sep=crcr]{%
19.875	31\\
20.125	31\\
};
\addplot [color=black, forget plot]
  table[row sep=crcr]{%
20.875	32\\
21.125	32\\
};
\addplot [color=black, forget plot]
  table[row sep=crcr]{%
21.875	38\\
22.125	38\\
};
\addplot [color=black, forget plot]
  table[row sep=crcr]{%
22.875	34\\
23.125	34\\
};
\addplot [color=black, forget plot]
  table[row sep=crcr]{%
23.875	38\\
24.125	38\\
};
\addplot [color=black, forget plot]
  table[row sep=crcr]{%
24.875	38\\
25.125	38\\
};
\addplot [color=black, forget plot]
  table[row sep=crcr]{%
25.875	42\\
26.125	42\\
};
\addplot [color=black, forget plot]
  table[row sep=crcr]{%
26.875	40\\
27.125	40\\
};
\addplot [color=black, forget plot]
  table[row sep=crcr]{%
27.875	43\\
28.125	43\\
};
\addplot [color=black, forget plot]
  table[row sep=crcr]{%
28.875	44\\
29.125	44\\
};
\addplot [color=black, forget plot]
  table[row sep=crcr]{%
29.875	43\\
30.125	43\\
};
\addplot [color=black, forget plot]
  table[row sep=crcr]{%
0.875	nan\\
1.125	nan\\
};
\addplot [color=black, forget plot]
  table[row sep=crcr]{%
1.875	5\\
2.125	5\\
};
\addplot [color=black, forget plot]
  table[row sep=crcr]{%
2.875	5\\
3.125	5\\
};
\addplot [color=black, forget plot]
  table[row sep=crcr]{%
3.875	6\\
4.125	6\\
};
\addplot [color=black, forget plot]
  table[row sep=crcr]{%
4.875	6\\
5.125	6\\
};
\addplot [color=black, forget plot]
  table[row sep=crcr]{%
5.875	8\\
6.125	8\\
};
\addplot [color=black, forget plot]
  table[row sep=crcr]{%
6.875	10\\
7.125	10\\
};
\addplot [color=black, forget plot]
  table[row sep=crcr]{%
7.875	11\\
8.125	11\\
};
\addplot [color=black, forget plot]
  table[row sep=crcr]{%
8.875	13\\
9.125	13\\
};
\addplot [color=black, forget plot]
  table[row sep=crcr]{%
9.875	15\\
10.125	15\\
};
\addplot [color=black, forget plot]
  table[row sep=crcr]{%
10.875	15\\
11.125	15\\
};
\addplot [color=black, forget plot]
  table[row sep=crcr]{%
11.875	16\\
12.125	16\\
};
\addplot [color=black, forget plot]
  table[row sep=crcr]{%
12.875	17\\
13.125	17\\
};
\addplot [color=black, forget plot]
  table[row sep=crcr]{%
13.875	18\\
14.125	18\\
};
\addplot [color=black, forget plot]
  table[row sep=crcr]{%
14.875	19\\
15.125	19\\
};
\addplot [color=black, forget plot]
  table[row sep=crcr]{%
15.875	20\\
16.125	20\\
};
\addplot [color=black, forget plot]
  table[row sep=crcr]{%
16.875	21\\
17.125	21\\
};
\addplot [color=black, forget plot]
  table[row sep=crcr]{%
17.875	22\\
18.125	22\\
};
\addplot [color=black, forget plot]
  table[row sep=crcr]{%
18.875	23\\
19.125	23\\
};
\addplot [color=black, forget plot]
  table[row sep=crcr]{%
19.875	23\\
20.125	23\\
};
\addplot [color=black, forget plot]
  table[row sep=crcr]{%
20.875	25\\
21.125	25\\
};
\addplot [color=black, forget plot]
  table[row sep=crcr]{%
21.875	25\\
22.125	25\\
};
\addplot [color=black, forget plot]
  table[row sep=crcr]{%
22.875	26\\
23.125	26\\
};
\addplot [color=black, forget plot]
  table[row sep=crcr]{%
23.875	26\\
24.125	26\\
};
\addplot [color=black, forget plot]
  table[row sep=crcr]{%
24.875	27\\
25.125	27\\
};
\addplot [color=black, forget plot]
  table[row sep=crcr]{%
25.875	28\\
26.125	28\\
};
\addplot [color=black, forget plot]
  table[row sep=crcr]{%
26.875	29\\
27.125	29\\
};
\addplot [color=black, forget plot]
  table[row sep=crcr]{%
27.875	29\\
28.125	29\\
};
\addplot [color=black, forget plot]
  table[row sep=crcr]{%
28.875	29\\
29.125	29\\
};
\addplot [color=black, forget plot]
  table[row sep=crcr]{%
29.875	31\\
30.125	31\\
};
\addplot [color=blue, forget plot]
  table[row sep=crcr]{%
0.75	nan\\
0.75	nan\\
1.25	nan\\
1.25	nan\\
0.75	nan\\
};
\addplot [color=blue, forget plot]
  table[row sep=crcr]{%
1.75	6\\
1.75	8\\
2.25	8\\
2.25	6\\
1.75	6\\
};
\addplot [color=blue, forget plot]
  table[row sep=crcr]{%
2.75	6\\
2.75	8\\
3.25	8\\
3.25	6\\
2.75	6\\
};
\addplot [color=blue, forget plot]
  table[row sep=crcr]{%
3.75	6\\
3.75	8.5\\
4.25	8.5\\
4.25	6\\
3.75	6\\
};
\addplot [color=blue, forget plot]
  table[row sep=crcr]{%
4.75	7\\
4.75	9\\
5.25	9\\
5.25	7\\
4.75	7\\
};
\addplot [color=blue, forget plot]
  table[row sep=crcr]{%
5.75	9\\
5.75	10\\
6.25	10\\
6.25	9\\
5.75	9\\
};
\addplot [color=blue, forget plot]
  table[row sep=crcr]{%
6.75	10\\
6.75	11\\
7.25	11\\
7.25	10\\
6.75	10\\
};
\addplot [color=blue, forget plot]
  table[row sep=crcr]{%
7.75	12\\
7.75	13\\
8.25	13\\
8.25	12\\
7.75	12\\
};
\addplot [color=blue, forget plot]
  table[row sep=crcr]{%
8.75	14\\
8.75	15\\
9.25	15\\
9.25	14\\
8.75	14\\
};
\addplot [color=blue, forget plot]
  table[row sep=crcr]{%
9.75	16\\
9.75	17\\
10.25	17\\
10.25	16\\
9.75	16\\
};
\addplot [color=blue, forget plot]
  table[row sep=crcr]{%
10.75	17\\
10.75	19\\
11.25	19\\
11.25	17\\
10.75	17\\
};
\addplot [color=blue, forget plot]
  table[row sep=crcr]{%
11.75	18\\
11.75	20\\
12.25	20\\
12.25	18\\
11.75	18\\
};
\addplot [color=blue, forget plot]
  table[row sep=crcr]{%
12.75	19.5\\
12.75	22\\
13.25	22\\
13.25	19.5\\
12.75	19.5\\
};
\addplot [color=blue, forget plot]
  table[row sep=crcr]{%
13.75	20.5\\
13.75	23\\
14.25	23\\
14.25	20.5\\
13.75	20.5\\
};
\addplot [color=blue, forget plot]
  table[row sep=crcr]{%
14.75	21.5\\
14.75	24\\
15.25	24\\
15.25	21.5\\
14.75	21.5\\
};
\addplot [color=blue, forget plot]
  table[row sep=crcr]{%
15.75	23\\
15.75	25\\
16.25	25\\
16.25	23\\
15.75	23\\
};
\addplot [color=blue, forget plot]
  table[row sep=crcr]{%
16.75	23.5\\
16.75	26\\
17.25	26\\
17.25	23.5\\
16.75	23.5\\
};
\addplot [color=blue, forget plot]
  table[row sep=crcr]{%
17.75	24.5\\
17.75	28\\
18.25	28\\
18.25	24.5\\
17.75	24.5\\
};
\addplot [color=blue, forget plot]
  table[row sep=crcr]{%
18.75	26\\
18.75	29\\
19.25	29\\
19.25	26\\
18.75	26\\
};
\addplot [color=blue, forget plot]
  table[row sep=crcr]{%
19.75	26.5\\
19.75	29\\
20.25	29\\
20.25	26.5\\
19.75	26.5\\
};
\addplot [color=blue, forget plot]
  table[row sep=crcr]{%
20.75	27\\
20.75	30\\
21.25	30\\
21.25	27\\
20.75	27\\
};
\addplot [color=blue, forget plot]
  table[row sep=crcr]{%
21.75	28\\
21.75	32\\
22.25	32\\
22.25	28\\
21.75	28\\
};
\addplot [color=blue, forget plot]
  table[row sep=crcr]{%
22.75	29\\
22.75	32\\
23.25	32\\
23.25	29\\
22.75	29\\
};
\addplot [color=blue, forget plot]
  table[row sep=crcr]{%
23.75	29.5\\
23.75	33\\
24.25	33\\
24.25	29.5\\
23.75	29.5\\
};
\addplot [color=blue, forget plot]
  table[row sep=crcr]{%
24.75	30.5\\
24.75	34.5\\
25.25	34.5\\
25.25	30.5\\
24.75	30.5\\
};
\addplot [color=blue, forget plot]
  table[row sep=crcr]{%
25.75	31\\
25.75	36\\
26.25	36\\
26.25	31\\
25.75	31\\
};
\addplot [color=blue, forget plot]
  table[row sep=crcr]{%
26.75	32\\
26.75	36\\
27.25	36\\
27.25	32\\
26.75	32\\
};
\addplot [color=blue, forget plot]
  table[row sep=crcr]{%
27.75	33\\
27.75	37\\
28.25	37\\
28.25	33\\
27.75	33\\
};
\addplot [color=blue, forget plot]
  table[row sep=crcr]{%
28.75	33.5\\
28.75	38\\
29.25	38\\
29.25	33.5\\
28.75	33.5\\
};
\addplot [color=blue, forget plot]
  table[row sep=crcr]{%
29.75	34\\
29.75	39\\
30.25	39\\
30.25	34\\
29.75	34\\
};
\addplot [color=red, forget plot]
  table[row sep=crcr]{%
0.75	nan\\
1.25	nan\\
};
\addplot [color=red, forget plot]
  table[row sep=crcr]{%
1.75	7\\
2.25	7\\
};
\addplot [color=red, forget plot]
  table[row sep=crcr]{%
2.75	7\\
3.25	7\\
};
\addplot [color=red, forget plot]
  table[row sep=crcr]{%
3.75	7\\
4.25	7\\
};
\addplot [color=red, forget plot]
  table[row sep=crcr]{%
4.75	8\\
5.25	8\\
};
\addplot [color=red, forget plot]
  table[row sep=crcr]{%
5.75	9\\
6.25	9\\
};
\addplot [color=red, forget plot]
  table[row sep=crcr]{%
6.75	11\\
7.25	11\\
};
\addplot [color=red, forget plot]
  table[row sep=crcr]{%
7.75	13\\
8.25	13\\
};
\addplot [color=red, forget plot]
  table[row sep=crcr]{%
8.75	15\\
9.25	15\\
};
\addplot [color=red, forget plot]
  table[row sep=crcr]{%
9.75	16\\
10.25	16\\
};
\addplot [color=red, forget plot]
  table[row sep=crcr]{%
10.75	18\\
11.25	18\\
};
\addplot [color=red, forget plot]
  table[row sep=crcr]{%
11.75	20\\
12.25	20\\
};
\addplot [color=red, forget plot]
  table[row sep=crcr]{%
12.75	21\\
13.25	21\\
};
\addplot [color=red, forget plot]
  table[row sep=crcr]{%
13.75	23\\
14.25	23\\
};
\addplot [color=red, forget plot]
  table[row sep=crcr]{%
14.75	24\\
15.25	24\\
};
\addplot [color=red, forget plot]
  table[row sep=crcr]{%
15.75	25\\
16.25	25\\
};
\addplot [color=red, forget plot]
  table[row sep=crcr]{%
16.75	26\\
17.25	26\\
};
\addplot [color=red, forget plot]
  table[row sep=crcr]{%
17.75	27\\
18.25	27\\
};
\addplot [color=red, forget plot]
  table[row sep=crcr]{%
18.75	28\\
19.25	28\\
};
\addplot [color=red, forget plot]
  table[row sep=crcr]{%
19.75	29\\
20.25	29\\
};
\addplot [color=red, forget plot]
  table[row sep=crcr]{%
20.75	30\\
21.25	30\\
};
\addplot [color=red, forget plot]
  table[row sep=crcr]{%
21.75	31\\
22.25	31\\
};
\addplot [color=red, forget plot]
  table[row sep=crcr]{%
22.75	32\\
23.25	32\\
};
\addplot [color=red, forget plot]
  table[row sep=crcr]{%
23.75	32\\
24.25	32\\
};
\addplot [color=red, forget plot]
  table[row sep=crcr]{%
24.75	34\\
25.25	34\\
};
\addplot [color=red, forget plot]
  table[row sep=crcr]{%
25.75	35\\
26.25	35\\
};
\addplot [color=red, forget plot]
  table[row sep=crcr]{%
26.75	36\\
27.25	36\\
};
\addplot [color=red, forget plot]
  table[row sep=crcr]{%
27.75	36\\
28.25	36\\
};
\addplot [color=red, forget plot]
  table[row sep=crcr]{%
28.75	37\\
29.25	37\\
};
\addplot [color=red, forget plot]
  table[row sep=crcr]{%
29.75	38\\
30.25	38\\
};
\addplot [color=black, only marks, mark=+, mark options={solid, draw=red}, forget plot]
  table[row sep=crcr]{%
nan	nan\\
};
\addplot [color=black, only marks, mark=+, mark options={solid, draw=red}, forget plot]
  table[row sep=crcr]{%
2	13\\
2	13\\
2	14\\
2	15\\
2	16\\
2	17\\
2	18\\
2	19\\
};
\addplot [color=black, only marks, mark=+, mark options={solid, draw=red}, forget plot]
  table[row sep=crcr]{%
3	12\\
3	13\\
3	15\\
3	16\\
3	17\\
3	18\\
3	20\\
3	21\\
3	22\\
};
\addplot [color=black, only marks, mark=+, mark options={solid, draw=red}, forget plot]
  table[row sep=crcr]{%
4	13\\
4	15\\
4	17\\
4	18\\
4	19\\
4	21\\
4	23\\
4	25\\
};
\addplot [color=black, only marks, mark=+, mark options={solid, draw=red}, forget plot]
  table[row sep=crcr]{%
5	15\\
5	17\\
5	19\\
5	21\\
5	23\\
5	26\\
5	27\\
};
\addplot [color=black, only marks, mark=+, mark options={solid, draw=red}, forget plot]
  table[row sep=crcr]{%
6	12\\
6	12\\
6	12\\
6	12\\
6	12\\
6	15\\
6	18\\
6	20\\
6	22\\
6	24\\
6	27\\
6	29\\
};
\addplot [color=black, only marks, mark=+, mark options={solid, draw=red}, forget plot]
  table[row sep=crcr]{%
7	13\\
7	13\\
7	13\\
7	13\\
7	14\\
7	18\\
7	21\\
7	24\\
7	26\\
7	29\\
7	31\\
};
\addplot [color=black, only marks, mark=+, mark options={solid, draw=red}, forget plot]
  table[row sep=crcr]{%
8	17\\
8	21\\
8	23\\
8	27\\
8	29\\
8	34\\
};
\addplot [color=black, only marks, mark=+, mark options={solid, draw=red}, forget plot]
  table[row sep=crcr]{%
9	17\\
9	17\\
9	22\\
9	25\\
9	28\\
9	31\\
9	35\\
};
\addplot [color=black, only marks, mark=+, mark options={solid, draw=red}, forget plot]
  table[row sep=crcr]{%
10	14\\
10	14\\
10	14\\
10	14\\
10	14\\
10	14\\
10	14\\
10	14\\
10	14\\
10	14\\
10	14\\
10	14\\
10	14\\
10	14\\
10	21\\
10	25\\
10	29\\
10	31\\
10	36\\
};
\addplot [color=black, only marks, mark=+, mark options={solid, draw=red}, forget plot]
  table[row sep=crcr]{%
11	26\\
11	29\\
11	32\\
11	37\\
};
\addplot [color=black, only marks, mark=+, mark options={solid, draw=red}, forget plot]
  table[row sep=crcr]{%
12	26\\
12	29\\
12	32\\
12	38\\
};
\addplot [color=black, only marks, mark=+, mark options={solid, draw=red}, forget plot]
  table[row sep=crcr]{%
13	27\\
13	29\\
13	33\\
13	40\\
};
\addplot [color=black, only marks, mark=+, mark options={solid, draw=red}, forget plot]
  table[row sep=crcr]{%
14	30\\
14	34\\
14	41\\
};
\addplot [color=black, only marks, mark=+, mark options={solid, draw=red}, forget plot]
  table[row sep=crcr]{%
15	30\\
15	34\\
15	42\\
};
\addplot [color=black, only marks, mark=+, mark options={solid, draw=red}, forget plot]
  table[row sep=crcr]{%
16	29\\
16	32\\
16	35\\
16	44\\
};
\addplot [color=black, only marks, mark=+, mark options={solid, draw=red}, forget plot]
  table[row sep=crcr]{%
17	30\\
17	32\\
17	36\\
17	45\\
};
\addplot [color=black, only marks, mark=+, mark options={solid, draw=red}, forget plot]
  table[row sep=crcr]{%
18	34\\
18	37\\
18	46\\
};
\addplot [color=black, only marks, mark=+, mark options={solid, draw=red}, forget plot]
  table[row sep=crcr]{%
19	35\\
19	37\\
19	47\\
};
\addplot [color=black, only marks, mark=+, mark options={solid, draw=red}, forget plot]
  table[row sep=crcr]{%
20	36\\
20	38\\
20	47\\
};
\addplot [color=black, only marks, mark=+, mark options={solid, draw=red}, forget plot]
  table[row sep=crcr]{%
21	35\\
21	39\\
21	48\\
};
\addplot [color=black, only marks, mark=+, mark options={solid, draw=red}, forget plot]
  table[row sep=crcr]{%
22	48\\
};
\addplot [color=black, only marks, mark=+, mark options={solid, draw=red}, forget plot]
  table[row sep=crcr]{%
23	37\\
23	39\\
23	49\\
};
\addplot [color=black, only marks, mark=+, mark options={solid, draw=red}, forget plot]
  table[row sep=crcr]{%
24	40\\
24	49\\
};
\addplot [color=black, only marks, mark=+, mark options={solid, draw=red}, forget plot]
  table[row sep=crcr]{%
25	41\\
25	50\\
};
\addplot [color=black, only marks, mark=+, mark options={solid, draw=red}, forget plot]
  table[row sep=crcr]{%
26	50\\
};
\addplot [color=black, only marks, mark=+, mark options={solid, draw=red}, forget plot]
  table[row sep=crcr]{%
27	43\\
27	50\\
};
\addplot [color=black, only marks, mark=+, mark options={solid, draw=red}, forget plot]
  table[row sep=crcr]{%
28	50\\
};
\addplot [color=black, only marks, mark=+, mark options={solid, draw=red}, forget plot]
  table[row sep=crcr]{%
29	50\\
};
\addplot [color=black, only marks, mark=+, mark options={solid, draw=red}, forget plot]
  table[row sep=crcr]{%
30	50\\
};
\end{axis}
\end{tikzpicture}%

%% file: iter_s5.tex
\definecolor{mycolor1}{rgb}{0.00000,0.44700,0.74100}%
\definecolor{mycolor2}{rgb}{0.85098,0.32549,0.09804}%
\begin{tikzpicture}

\begin{axis}[%
width=0.8\columnwidth,
height=0.8\columnwidth,
at={(0\columnwidth,0\columnwidth)},
scale only axis,
xmin=0,
xmax=100,
separate axis lines,
every outer y axis line/.append style={mycolor1},
every y tick label/.append style={font=\color{mycolor1}},
every y tick/.append style={mycolor1},
ymin=5,
ymax=7,
ylabel style={font=\color{mycolor1}\footnotesize},
ylabel={Newton Iteration},
axis background/.style={fill=none},
yticklabel style = {font=\color{mycolor1}\small},
yticklabel pos=left
]
\addplot [color=mycolor1, line width=2.0pt, mark=o, mark options={solid, mycolor1}, forget plot]
  table[row sep=crcr]{%
1	6\\
2	6\\
3	6\\
4	6\\
5	6\\
6	6\\
7	6\\
8	6\\
9	6\\
10	6\\
11	6\\
12	6\\
13	6\\
14	6\\
15	6\\
16	6\\
17	6\\
18	6\\
19	6\\
20	6\\
21	6\\
22	6\\
23	6\\
24	6\\
25	6\\
26	6\\
27	6\\
28	6\\
29	6\\
30	6\\
31	6\\
32	6\\
33	6\\
34	6\\
35	6\\
36	6\\
37	6\\
38	6\\
39	6\\
40	6\\
41	6\\
42	6\\
43	6\\
44	6\\
45	6\\
46	6\\
47	6\\
48	6\\
49	6\\
50	6\\
51	6\\
52	6\\
53	6\\
54	6\\
55	6\\
56	6\\
57	6\\
58	6\\
59	6\\
60	6\\
61	6\\
62	6\\
63	6\\
64	6\\
65	6\\
66	6\\
67	6\\
68	6\\
69	6\\
70	6\\
71	6\\
72	6\\
73	6\\
74	6\\
75	6\\
76	6\\
77	6\\
78	6\\
79	6\\
80	6\\
81	6\\
82	6\\
83	6\\
84	6\\
85	6\\
86	6\\
87	6\\
88	6\\
89	6\\
90	6\\
91	6\\
92	6\\
93	6\\
94	6\\
95	6\\
96	6\\
97	6\\
98	6\\
99	6\\
100	6\\
};
\end{axis}
\begin{axis}[%
width=0.8\columnwidth,
height=0.8\columnwidth,
at={(0\columnwidth,0\columnwidth)},
scale only axis,
xmin=0,
xmax=100,
separate axis lines,
every outer y axis line/.append style={mycolor2},
every y tick label/.append style={font=\color{mycolor2}},
every y tick/.append style={mycolor2},
yticklabels={},
ymin=10,
ymax=20,
ylabel style={font=\color{mycolor2}\footnotesize},
yticklabel pos=right
]
\addplot [color=mycolor2, line width=2.0pt, mark=x, mark options={solid, mycolor2}, forget plot]
  table[row sep=crcr]{%
1	17\\
2	17\\
3	14\\
4	14\\
5	14\\
6	14\\
7	14\\
8	14\\
9	14\\
10	14\\
11	14\\
12	14\\
13	14\\
14	14\\
15	14\\
16	14\\
17	14\\
18	14\\
19	14\\
20	14\\
21	14\\
22	14\\
23	14\\
24	14\\
25	14\\
26	14\\
27	14\\
28	14\\
29	14\\
30	14\\
31	14\\
32	14\\
33	14\\
34	14\\
35	14\\
36	14\\
37	14\\
38	14\\
39	14\\
40	14\\
41	14\\
42	12\\
43	12\\
44	12\\
45	12\\
46	12\\
47	12\\
48	12\\
49	12\\
50	12\\
51	12\\
52	12\\
53	12\\
54	12\\
55	12\\
56	12\\
57	12\\
58	12\\
59	12\\
60	12\\
61	12\\
62	12\\
63	12\\
64	12\\
65	12\\
66	12\\
67	12\\
68	12\\
69	12\\
70	12\\
71	12\\
72	12\\
73	12\\
74	12\\
75	12\\
76	12\\
77	12\\
78	12\\
79	12\\
80	12\\
81	12\\
82	12\\
83	12\\
84	12\\
85	12\\
86	12\\
87	12\\
88	12\\
89	12\\
90	12\\
91	12\\
92	12\\
93	12\\
94	12\\
95	12\\
96	12\\
97	12\\
98	12\\
99	12\\
100	12\\
};
\end{axis}
\end{tikzpicture}%

%% file: iter_s15.tex
\definecolor{mycolor1}{rgb}{0.00000,0.44700,0.74100}%
\definecolor{mycolor2}{rgb}{0.85098,0.32549,0.09804}%
\begin{tikzpicture}

\begin{axis}[%
width=0.8\columnwidth,
height=0.8\columnwidth,
at={(0\columnwidth,0\columnwidth)},
scale only axis,
xmin=0,
xmax=100,
separate axis lines,
every outer y axis line/.append style={mycolor1},
every y tick label/.append style={font=\color{mycolor1}},
every y tick/.append style={mycolor1},
ymin=5,
ymax=7,
ylabel style={font=\color{mycolor1}\footnotesize},
axis background/.style={fill=none},
yticklabels={},
yticklabel pos=left
]
\addplot [color=mycolor1, line width=2.0pt, mark=o, mark options={solid, mycolor1}, forget plot]
  table[row sep=crcr]{%
1	6\\
2	6\\
3	6\\
4	6\\
5	6\\
6	6\\
7	6\\
8	6\\
9	6\\
10	6\\
11	6\\
12	6\\
13	6\\
14	6\\
15	6\\
16	6\\
17	6\\
18	6\\
19	6\\
20	6\\
21	6\\
22	6\\
23	6\\
24	6\\
25	6\\
26	6\\
27	6\\
28	6\\
29	6\\
30	6\\
31	6\\
32	6\\
33	6\\
34	6\\
35	6\\
36	6\\
37	6\\
38	6\\
39	6\\
40	6\\
41	6\\
42	6\\
43	6\\
44	6\\
45	6\\
46	6\\
47	6\\
48	6\\
49	6\\
50	6\\
51	6\\
52	6\\
53	6\\
54	6\\
55	6\\
56	6\\
57	6\\
58	6\\
59	6\\
60	6\\
61	6\\
62	6\\
63	6\\
64	6\\
65	6\\
66	6\\
67	6\\
68	6\\
69	6\\
70	6\\
71	6\\
72	6\\
73	6\\
74	6\\
75	6\\
76	6\\
77	6\\
78	6\\
79	6\\
80	6\\
81	6\\
82	6\\
83	6\\
84	6\\
85	6\\
86	6\\
87	6\\
88	6\\
89	6\\
90	6\\
91	6\\
92	6\\
93	6\\
94	6\\
95	6\\
96	6\\
97	6\\
98	6\\
99	6\\
100	6\\
};
\end{axis}
\begin{axis}[%
width=0.8\columnwidth,
height=0.8\columnwidth,
at={(0\columnwidth,0\columnwidth)},
scale only axis,
xmin=0,
xmax=100,
separate axis lines,
every outer y axis line/.append style={mycolor2},
every y tick label/.append style={font=\color{mycolor2}},
every y tick/.append style={mycolor2},
ymin=10,
ymax=20,
ylabel style={font=\color{mycolor2}\footnotesize},
yticklabels={},
yticklabel pos=right
]
\addplot [color=mycolor2, line width=2.0pt, mark=x, mark options={solid, mycolor2}, forget plot]
  table[row sep=crcr]{%
1	18\\
2	15\\
3	14\\
4	14\\
5	14\\
6	14\\
7	14\\
8	14\\
9	14\\
10	14\\
11	14\\
12	14\\
13	14\\
14	14\\
15	14\\
16	14\\
17	14\\
18	14\\
19	14\\
20	14\\
21	14\\
22	14\\
23	14\\
24	14\\
25	14\\
26	14\\
27	14\\
28	14\\
29	14\\
30	14\\
31	14\\
32	14\\
33	14\\
34	14\\
35	14\\
36	14\\
37	14\\
38	14\\
39	14\\
40	14\\
41	14\\
42	14\\
43	14\\
44	14\\
45	14\\
46	14\\
47	14\\
48	14\\
49	14\\
50	14\\
51	14\\
52	14\\
53	14\\
54	14\\
55	14\\
56	14\\
57	14\\
58	14\\
59	14\\
60	14\\
61	14\\
62	14\\
63	14\\
64	14\\
65	14\\
66	14\\
67	14\\
68	14\\
69	14\\
70	14\\
71	14\\
72	14\\
73	14\\
74	14\\
75	14\\
76	14\\
77	14\\
78	14\\
79	14\\
80	14\\
81	14\\
82	14\\
83	14\\
84	14\\
85	14\\
86	14\\
87	14\\
88	14\\
89	14\\
90	14\\
91	14\\
92	13\\
93	13\\
94	13\\
95	13\\
96	13\\
97	13\\
98	13\\
99	13\\
100	13\\
};
\end{axis}
\end{tikzpicture}%

%% file: iter_s25.tex
\definecolor{mycolor1}{rgb}{0.00000,0.44700,0.74100}%
\definecolor{mycolor2}{rgb}{0.85098,0.32549,0.09804}%
\begin{tikzpicture}

\begin{axis}[%
width=0.8\columnwidth,
height=0.8\columnwidth,
at={(0\columnwidth,0\columnwidth)},
scale only axis,
xmin=0,
xmax=100,
separate axis lines,
every outer y axis line/.append style={mycolor1},
every y tick label/.append style={font=\color{mycolor1}},
every y tick/.append style={mycolor1},
ymin=5,
ymax=7,
ylabel style={font=\color{mycolor1}\footnotesize},
axis background/.style={fill=none},
yticklabels={},
yticklabel pos=left
]
\addplot [color=mycolor1, line width=2.0pt, mark=o, mark options={solid, mycolor1}, forget plot]
  table[row sep=crcr]{%
1	6\\
2	6\\
3	6\\
4	6\\
5	6\\
6	6\\
7	6\\
8	6\\
9	6\\
10	6\\
11	6\\
12	6\\
13	6\\
14	6\\
15	6\\
16	6\\
17	6\\
18	6\\
19	6\\
20	6\\
21	6\\
22	6\\
23	6\\
24	6\\
25	6\\
26	6\\
27	6\\
28	6\\
29	6\\
30	6\\
31	6\\
32	6\\
33	6\\
34	6\\
35	6\\
36	6\\
37	6\\
38	6\\
39	6\\
40	6\\
41	6\\
42	6\\
43	6\\
44	6\\
45	6\\
46	6\\
47	6\\
48	6\\
49	6\\
50	6\\
51	6\\
52	6\\
53	6\\
54	6\\
55	6\\
56	6\\
57	6\\
58	6\\
59	6\\
60	6\\
61	6\\
62	6\\
63	6\\
64	6\\
65	6\\
66	6\\
67	6\\
68	6\\
69	6\\
70	6\\
71	6\\
72	6\\
73	6\\
74	6\\
75	6\\
76	6\\
77	6\\
78	6\\
79	6\\
80	6\\
81	6\\
82	6\\
83	6\\
84	6\\
85	6\\
86	6\\
87	6\\
88	6\\
89	6\\
90	6\\
91	6\\
92	6\\
93	6\\
94	6\\
95	6\\
96	6\\
97	6\\
98	6\\
99	6\\
100	6\\
};
\end{axis}
\begin{axis}[%
width=0.8\columnwidth,
height=0.8\columnwidth,
at={(0\columnwidth,0\columnwidth)},
scale only axis,
xmin=0,
xmax=100,
separate axis lines,
every outer y axis line/.append style={mycolor2},
every y tick label/.append style={font=\color{mycolor2}},
every y tick/.append style={mycolor2},
ymin=10,
ymax=20,
ylabel style={font=\color{mycolor2}\footnotesize},
ylabel={Krylov Dimension},
yticklabel style = {font=\color{mycolor2}\small},
yticklabel pos=right
]
\addplot [color=mycolor2, line width=2.0pt, mark=x, mark options={solid, mycolor2}, forget plot]
  table[row sep=crcr]{%
1	18\\
2	17\\
3	14\\
4	14\\
5	14\\
6	14\\
7	14\\
8	14\\
9	14\\
10	14\\
11	14\\
12	14\\
13	14\\
14	14\\
15	14\\
16	14\\
17	14\\
18	14\\
19	14\\
20	14\\
21	14\\
22	14\\
23	14\\
24	14\\
25	14\\
26	14\\
27	14\\
28	14\\
29	14\\
30	14\\
31	14\\
32	14\\
33	14\\
34	14\\
35	14\\
36	14\\
37	14\\
38	14\\
39	14\\
40	14\\
41	14\\
42	14\\
43	14\\
44	14\\
45	14\\
46	14\\
47	14\\
48	14\\
49	14\\
50	14\\
51	14\\
52	14\\
53	14\\
54	14\\
55	14\\
56	14\\
57	14\\
58	14\\
59	14\\
60	14\\
61	14\\
62	14\\
63	14\\
64	14\\
65	14\\
66	14\\
67	14\\
68	14\\
69	14\\
70	14\\
71	14\\
72	14\\
73	14\\
74	14\\
75	14\\
76	14\\
77	14\\
78	14\\
79	14\\
80	14\\
81	14\\
82	14\\
83	14\\
84	14\\
85	14\\
86	14\\
87	14\\
88	14\\
89	14\\
90	14\\
91	14\\
92	14\\
93	14\\
94	14\\
95	14\\
96	14\\
97	13\\
98	13\\
99	13\\
100	13\\
};
\end{axis}
\end{tikzpicture}%

%% file: nonlinearperformance.tex
\definecolor{mycolor1}{rgb}{0.00000,0.44706,0.74118}%
\definecolor{mycolor2}{rgb}{0.85098,0.32549,0.09804}%
\begin{tikzpicture}

\begin{axis}[%
width=0.788\columnwidth,
height=0.196\columnwidth,
at={(0\columnwidth,0\columnwidth)},
scale only axis,
xmin=0,
xmax=30,
xlabel style={font=\color{white!15!black}},
xlabel={s : number of stages},
separate axis lines,
every outer y axis line/.append style={mycolor1},
every y tick label/.append style={font=\color{mycolor1}},
every y tick/.append style={mycolor1},
ymin=0,
ymax=60,
ylabel style={font=\color{mycolor1}},
ylabel={Elapsed Time (s)},
axis background/.style={fill=none},
yticklabel pos=left,
xmajorgrids,
ymajorgrids,
legend style={at={(0.03,1.01)}, anchor=north west, legend cell align=left, align=left, draw=none, fill=none}
]
\addplot [color=mycolor1, line width=2.0pt, mark=square, mark options={solid, mycolor1}]
  table[row sep=crcr]{%
1	1.283\\
2	3.583\\
3	5.302\\
4	4.63\\
5	4.851\\
6	5.787\\
7	6.121\\
8	7.03\\
9	8.133\\
10	7.473\\
11	9.368\\
12	11.261\\
13	8.317\\
14	9.482\\
15	11.604\\
16	9.196\\
17	12.087\\
18	10.938\\
19	13.501\\
20	16.732\\
21	14.313\\
22	15.311\\
23	13.57\\
24	18.985\\
25	21.103\\
26	18.535\\
27	19.45\\
28	17.33\\
29	17.465\\
30	18.121\\
};
\addlegendentry{$s$ threads}
\addplot [color=mycolor1, line width=2.0pt, mark=o, mark options={solid, mycolor1}]
  table[row sep=crcr]{%
1	1.283\\
2	3.142\\
3	3.482\\
4	6.057\\
5	6.786\\
6	9.226\\
7	9.777\\
8	12.444\\
9	12.895\\
10	15.804\\
11	16.454\\
12	19.379\\
13	19.788\\
14	22.985\\
15	24.113\\
16	26.184\\
17	27.886\\
18	29.897\\
19	32.587\\
20	34.455\\
21	35.919\\
22	38.363\\
23	40.808\\
24	42.134\\
25	44.607\\
26	46.939\\
27	50.751\\
28	52.59\\
29	55.841\\
30	56.921\\
};
\addlegendentry{sequential ($1$ thread)}
\end{axis}
\begin{axis}[%
width=0.788\columnwidth,
height=0.196\columnwidth,
at={(0\columnwidth,0\columnwidth)},
scale only axis,
xmin=0,
xmax=30,
xlabel style={font=\color{white!15!black}},
separate axis lines,
every outer y axis line/.append style={mycolor2},
every y tick label/.append style={font=\color{mycolor2}},
every y tick/.append style={mycolor2},
ymin=0.5,
ymax=3.5,
ylabel style={font=\color{mycolor2}},
ylabel={Speedup},
axis background/.style={fill=none},
yticklabel pos=right,
xmajorgrids,
ymajorgrids,
legend style={at={(0.03,0.68)}, anchor=north west, legend cell align=left, align=left, draw=none, fill=none}
]
\addplot [color=mycolor2, line width=2.0pt] %
  table[row sep=crcr]{%
1	1\\
2	0.876918783142618\\
3	0.65673330818559\\
4	1.30820734341253\\
5	1.39888682745826\\
6	1.59426300328322\\
7	1.59728802483254\\
8	1.7701280227596\\
9	1.58551579982786\\
10	2.11481332798073\\
11	1.7564047822374\\
12	1.72089512476689\\
13	2.37922327762414\\
14	2.42406665260494\\
15	2.07799034815581\\
16	2.84732492387995\\
17	2.30710680896831\\
18	2.73331504845493\\
19	2.41367306125472\\
20	2.05922782691848\\
21	2.50953678474114\\
22	2.50558422049507\\
23	3.0072218128224\\
24	2.2193310508296\\
25	2.11377529261242\\
26	2.53245211761532\\
27	2.6093059125964\\
28	3.03462204270052\\
29	3.19730890352133\\
30	3.14116218751725\\
};
\addlegendentry{Speedup}

\end{axis}
\end{tikzpicture}%

%% file: main-revision.bbl

\begin{thebibliography}{34}
\ifx \bisbn   \undefined \def \bisbn  #1{ISBN #1}\fi
\ifx \binits  \undefined \def \binits#1{#1}\fi
\ifx \bauthor  \undefined \def \bauthor#1{#1}\fi
\ifx \batitle  \undefined \def \batitle#1{#1}\fi
\ifx \bjtitle  \undefined \def \bjtitle#1{#1}\fi
\ifx \bvolume  \undefined \def \bvolume#1{\textbf{#1}}\fi
\ifx \byear  \undefined \def \byear#1{#1}\fi
\ifx \bissue  \undefined \def \bissue#1{#1}\fi
\ifx \bfpage  \undefined \def \bfpage#1{#1}\fi
\ifx \blpage  \undefined \def \blpage #1{#1}\fi
\ifx \burl  \undefined \def \burl#1{\textsf{#1}}\fi
\ifx \doiurl  \undefined \def \doiurl#1{\url{https://doi.org/#1}}\fi
\ifx \betal  \undefined \def \betal{\textit{et al.}}\fi
\ifx \binstitute  \undefined \def \binstitute#1{#1}\fi
\ifx \binstitutionaled  \undefined \def \binstitutionaled#1{#1}\fi
\ifx \bctitle  \undefined \def \bctitle#1{#1}\fi
\ifx \beditor  \undefined \def \beditor#1{#1}\fi
\ifx \bpublisher  \undefined \def \bpublisher#1{#1}\fi
\ifx \bbtitle  \undefined \def \bbtitle#1{#1}\fi
\ifx \bedition  \undefined \def \bedition#1{#1}\fi
\ifx \bseriesno  \undefined \def \bseriesno#1{#1}\fi
\ifx \blocation  \undefined \def \blocation#1{#1}\fi
\ifx \bsertitle  \undefined \def \bsertitle#1{#1}\fi
\ifx \bsnm \undefined \def \bsnm#1{#1}\fi
\ifx \bsuffix \undefined \def \bsuffix#1{#1}\fi
\ifx \bparticle \undefined \def \bparticle#1{#1}\fi
\ifx \barticle \undefined \def \barticle#1{#1}\fi
\bibcommenthead
\ifx \bconfdate \undefined \def \bconfdate #1{#1}\fi
\ifx \botherref \undefined \def \botherref #1{#1}\fi
\ifx \url \undefined \def \url#1{\textsf{#1}}\fi
\ifx \bchapter \undefined \def \bchapter#1{#1}\fi
\ifx \bbook \undefined \def \bbook#1{#1}\fi
\ifx \bcomment \undefined \def \bcomment#1{#1}\fi
\ifx \oauthor \undefined \def \oauthor#1{#1}\fi
\ifx \citeauthoryear \undefined \def \citeauthoryear#1{#1}\fi
\ifx \endbibitem  \undefined \def \endbibitem {}\fi
\ifx \bconflocation  \undefined \def \bconflocation#1{#1}\fi
\ifx \arxivurl  \undefined \def \arxivurl#1{\textsf{#1}}\fi
\csname PreBibitemsHook\endcsname

\bibitem[\protect\citeauthoryear{Badia and Verdugo}{2020}]{Badia2020}
\begin{barticle}
\bauthor{\bsnm{Badia}, \binits{S.}},
\bauthor{\bsnm{Verdugo}, \binits{F.}}:
\batitle{Gridap: {A}n extensible {F}inite {E}lement toolbox in {J}ulia}.
\bjtitle{Journal of Open Source Software}
\bvolume{5}(\bissue{52}),
\bfpage{2520}
(\byear{2020})
\doiurl{10.21105/joss.02520}
\end{barticle}
\endbibitem

\bibitem[\protect\citeauthoryear{Bellen}{1989}]{Bellen}
\begin{bchapter}
\bauthor{\bsnm{Bellen}, \binits{A.}}:
\bctitle{Parallelism across the steps for difference and differential
  equations}.
In: \bbtitle{Numerical Methods for Ordinary Differential Equations
  ({L}'{A}quila, 1987)}.
\bsertitle{Lecture Notes in Math.},
vol. \bseriesno{1386},
pp. \bfpage{22}--\blpage{35}.
\bpublisher{Springer},
\blocation{Berlin}
(\byear{1989}).
\doiurl{10.1007/BFb0089229} .
\burl{https://doi.org/10.1007/BFb0089229}
\end{bchapter}
\endbibitem

\bibitem[\protect\citeauthoryear{Bellen et~al.}{1990}]{Bellen2}
\begin{barticle}
\bauthor{\bsnm{Bellen}, \binits{A.}},
\bauthor{\bsnm{Vermiglio}, \binits{R.}},
\bauthor{\bsnm{Zennaro}, \binits{M.}}:
\batitle{Parallel {ODE}-solvers with stepsize control}.
\bjtitle{J. Comput. Appl. Math.}
\bvolume{31}(\bissue{2}),
\bfpage{277}--\blpage{293}
(\byear{1990})
\doiurl{10.1016/0377-0427(90)90170-5}
\end{barticle}
\endbibitem

\bibitem[\protect\citeauthoryear{Berljafa and G\"uttel}{2015}]{MR3361437}
\begin{barticle}
\bauthor{\bsnm{Berljafa}, \binits{M.}},
\bauthor{\bsnm{G\"uttel}, \binits{S.}}:
\batitle{Generalized rational {K}rylov decompositions with an application to
  rational approximation}.
\bjtitle{SIAM J. Matrix Anal. Appl.}
\bvolume{36}(\bissue{2}),
\bfpage{894}--\blpage{916}
(\byear{2015})
\doiurl{10.1137/140998081}
\end{barticle}
\endbibitem

\bibitem[\protect\citeauthoryear{Brugnano and Iavernaro}{2016}]{BrugnanoBook}
\begin{bbook}
\bauthor{\bsnm{Brugnano}, \binits{L.}},
\bauthor{\bsnm{Iavernaro}, \binits{F.}}:
\bbtitle{Line Integral Methods for Conservative Problems}.
\bsertitle{Monographs and Research Notes in Mathematics},
p. \bfpage{222}.
\bpublisher{CRC Press},
\blocation{Boca Raton, FL}
(\byear{2016}).
\doiurl{10.1201/b19319}
\end{bbook}
\endbibitem

\bibitem[\protect\citeauthoryear{Collar}{1962}]{CENTROSKEW}
\begin{barticle}
\bauthor{\bsnm{Collar}, \binits{A.R.}}:
\batitle{On centrosymmetric and centroskew matrices}.
\bjtitle{Quart. J. Mech. Appl. Math.}
\bvolume{15},
\bfpage{265}--\blpage{281}
(\byear{1962})
\doiurl{10.1093/qjmam/15.3.265}
\end{barticle}
\endbibitem

\bibitem[\protect\citeauthoryear{Cong}{1993}]{Cong1}
\begin{barticle}
\bauthor{\bsnm{Cong}, \binits{N.H.}}:
\batitle{Note on the performance of direct and indirect
  {R}unge-{K}utta-{N}ystr\"{o}m methods}.
\bjtitle{J. Comput. Appl. Math.}
\bvolume{45}(\bissue{3}),
\bfpage{347}--\blpage{355}
(\byear{1993})
\doiurl{10.1016/0377-0427(93)90053-E}
\end{barticle}
\endbibitem

\bibitem[\protect\citeauthoryear{Cong et~al.}{1999}]{Cong2}
\begin{barticle}
\bauthor{\bsnm{Cong}, \binits{N.H.}},
\bauthor{\bsnm{Strehmel}, \binits{K.}},
\bauthor{\bsnm{Weiner}, \binits{R.}},
\bauthor{\bsnm{Podhaisky}, \binits{H.}}:
\batitle{Runge-{K}utta-{N}ystr\"{o}m-type parallel block predictor-corrector
  methods}.
\bjtitle{Adv. Comput. Math.}
\bvolume{10}(\bissue{2}),
\bfpage{115}--\blpage{133}
(\byear{1999})
\doiurl{10.1023/A:1018930732643}
\end{barticle}
\endbibitem

\bibitem[\protect\citeauthoryear{D'Ambra et~al.}{2023}]{DAmbra2023}
\begin{botherref}
\oauthor{\bsnm{D'Ambra}, \binits{P.}},
\oauthor{\bsnm{Durastante}, \binits{F.}},
\oauthor{\bsnm{Filippone}, \binits{S.}}:
Parallel sparse computation toolkit[formula presented].
Software Impacts
\textbf{15}
(2023)
\doiurl{10.1016/j.simpa.2022.100463} .
Cited by: 1; All Open Access, Gold Open Access, Green Open Access
\end{botherref}
\endbibitem

\bibitem[\protect\citeauthoryear{Davis}{2004}]{MR2075981}
\begin{barticle}
\bauthor{\bsnm{Davis}, \binits{T.A.}}:
\batitle{Algorithm 832: {UMFPACK} {V}4.3---an unsymmetric-pattern multifrontal
  method}.
\bjtitle{ACM Trans. Math. Software}
\bvolume{30}(\bissue{2}),
\bfpage{196}--\blpage{199}
(\byear{2004})
\doiurl{10.1145/992200.992206}
\end{barticle}
\endbibitem

\bibitem[\protect\citeauthoryear{Gander and Palitta}{2024}]{MR4713232}
\begin{barticle}
\bauthor{\bsnm{Gander}, \binits{M.J.}},
\bauthor{\bsnm{Palitta}, \binits{D.}}:
\batitle{A new {P}ara{D}iag time-parallel time integration method}.
\bjtitle{SIAM J. Sci. Comput.}
\bvolume{46}(\bissue{2}),
\bfpage{697}--\blpage{718}
(\byear{2024})
\doiurl{10.1137/23M1568028}
\end{barticle}
\endbibitem

\bibitem[\protect\citeauthoryear{Gander and Wu}{2020}]{MR4167091}
\begin{barticle}
\bauthor{\bsnm{Gander}, \binits{M.J.}},
\bauthor{\bsnm{Wu}, \binits{S.-L.}}:
\batitle{A diagonalization-based parareal algorithm for dissipative and wave
  propagation problems}.
\bjtitle{SIAM J. Numer. Anal.}
\bvolume{58}(\bissue{5}),
\bfpage{2981}--\blpage{3009}
(\byear{2020})
\doiurl{10.1137/19M1271683}
\end{barticle}
\endbibitem

\bibitem[\protect\citeauthoryear{Gautschi}{2004}]{Gautschi}
\begin{bbook}
\bauthor{\bsnm{Gautschi}, \binits{W.}}:
\bbtitle{Orthogonal Polynomials: Computation and Approximation}.
\bsertitle{Numerical Mathematics and Scientific Computation},
p. \bfpage{301}.
\bpublisher{Oxford University Press},
\blocation{New York}
(\byear{2004}).
\doiurl{10.1093/oso/9780198506720.001.0001}
\end{bbook}
\endbibitem

\bibitem[\protect\citeauthoryear{Hairer and Wanner}{1996}]{BibleVolII}
\begin{bbook}
\bauthor{\bsnm{Hairer}, \binits{E.}},
\bauthor{\bsnm{Wanner}, \binits{G.}}:
\bbtitle{Solving Ordinary Differential Equations. {II}},
\bedition{2}nd edn.
\bsertitle{Springer Series in Computational Mathematics},
vol. \bseriesno{14},
p. \bfpage{614}.
\bpublisher{Springer},
\blocation{Berlin}
(\byear{1996}).
\doiurl{10.1007/978-3-642-05221-7}
\end{bbook}
\endbibitem

\bibitem[\protect\citeauthoryear{Hairer et~al.}{1993}]{BibleVol1}
\begin{bbook}
\bauthor{\bsnm{Hairer}, \binits{E.}},
\bauthor{\bsnm{N{\o}rsett}, \binits{S.P.}},
\bauthor{\bsnm{Wanner}, \binits{G.}}:
\bbtitle{Solving Ordinary Differential Equations. {I}},
\bedition{2}nd edn.
\bsertitle{Springer Series in Computational Mathematics},
vol. \bseriesno{8},
p. \bfpage{528}.
\bpublisher{Springer},
\blocation{Berlin}
(\byear{1993}).
\doiurl{10.1007/978-3-540-78862-1}
\end{bbook}
\endbibitem

\bibitem[\protect\citeauthoryear{Hairer et~al.}{2002}]{BibleVolIII}
\begin{bbook}
\bauthor{\bsnm{Hairer}, \binits{E.}},
\bauthor{\bsnm{Lubich}, \binits{C.}},
\bauthor{\bsnm{Wanner}, \binits{G.}}:
\bbtitle{Geometric Numerical Integration}.
\bsertitle{Springer Series in Computational Mathematics},
vol. \bseriesno{31},
p. \bfpage{515}.
\bpublisher{Springer},
\blocation{Berlin}
(\byear{2002}).
\doiurl{10.1007/978-3-662-05018-7}
\end{bbook}
\endbibitem

\bibitem[\protect\citeauthoryear{Iserles and N{\o}rsett}{1990}]{Iserles}
\begin{barticle}
\bauthor{\bsnm{Iserles}, \binits{A.}},
\bauthor{\bsnm{N{\o}rsett}, \binits{S.P.}}:
\batitle{On the theory of parallel {R}unge-{K}utta methods}.
\bjtitle{IMA J. Numer. Anal.}
\bvolume{10}(\bissue{4}),
\bfpage{463}--\blpage{488}
(\byear{1990})
\doiurl{10.1093/imanum/10.4.463}
\end{barticle}
\endbibitem

\bibitem[\protect\citeauthoryear{Jay}{2000}]{MR1790038}
\begin{barticle}
\bauthor{\bsnm{Jay}, \binits{L.O.}}:
\batitle{Inexact simplified {N}ewton iterations for implicit {R}unge-{K}utta
  methods}.
\bjtitle{SIAM J. Numer. Anal.}
\bvolume{38}(\bissue{4}),
\bfpage{1369}--\blpage{1388}
(\byear{2000})
\doiurl{10.1137/S0036142999360573}
\end{barticle}
\endbibitem

\bibitem[\protect\citeauthoryear{K{\aa}gstr\"om and Poromaa}{1996}]{MR1383186}
\begin{barticle}
\bauthor{\bsnm{K{\aa}gstr\"om}, \binits{B.}},
\bauthor{\bsnm{Poromaa}, \binits{P.}}:
\batitle{L{APACK}-style algorithms and software for solving the generalized
  {S}ylvester equation and estimating the separation between regular matrix
  pairs}.
\bjtitle{ACM Trans. Math. Software}
\bvolume{22}(\bissue{1}),
\bfpage{78}--\blpage{103}
(\byear{1996})
\doiurl{10.1145/225545.225552}
\end{barticle}
\endbibitem

\bibitem[\protect\citeauthoryear{Kressner et~al.}{2023}]{MR4646959}
\begin{barticle}
\bauthor{\bsnm{Kressner}, \binits{D.}},
\bauthor{\bsnm{Massei}, \binits{S.}},
\bauthor{\bsnm{Zhu}, \binits{J.}}:
\batitle{Improved {P}ara{D}iag via low-rank updates and interpolation}.
\bjtitle{Numer. Math.}
\bvolume{155}(\bissue{1-2}),
\bfpage{175}--\blpage{209}
(\byear{2023})
\doiurl{10.1007/s00211-023-01372-w}
\end{barticle}
\endbibitem

\bibitem[\protect\citeauthoryear{Leveque et~al.}{2024}]{SantoloSVD}
\begin{barticle}
\bauthor{\bsnm{Leveque}, \binits{S.}},
\bauthor{\bsnm{Bergamaschi}, \binits{L.}},
\bauthor{\bsnm{Mart\'inez}, \binits{A.}},
\bauthor{\bsnm{Pearson}, \binits{J.W.}}:
\batitle{Parallel-in-time solver for the all-at-once {R}unge-{K}utta
  discretization}.
\bjtitle{SIAM J. Matrix Anal. Appl.}
\bvolume{45}(\bissue{4}),
\bfpage{1902}--\blpage{1928}
(\byear{2024})
\doiurl{10.1137/23M1567862}
\end{barticle}
\endbibitem

\bibitem[\protect\citeauthoryear{Liniger and
  Willoughby}{1970}]{SimplifiedNewton}
\begin{barticle}
\bauthor{\bsnm{Liniger}, \binits{W.}},
\bauthor{\bsnm{Willoughby}, \binits{R.A.}}:
\batitle{Efficient integration methods for stiff systems of ordinary
  differential equations}.
\bjtitle{SIAM J. Numer. Anal.}
\bvolume{7},
\bfpage{47}--\blpage{66}
(\byear{1970})
\doiurl{10.1137/0707002}
\end{barticle}
\endbibitem

\bibitem[\protect\citeauthoryear{Munch et~al.}{2024}]{MR4735250}
\begin{barticle}
\bauthor{\bsnm{Munch}, \binits{P.}},
\bauthor{\bsnm{Dravins}, \binits{I.}},
\bauthor{\bsnm{Kronbichler}, \binits{M.}},
\bauthor{\bsnm{Neytcheva}, \binits{M.}}:
\batitle{Stage-parallel fully implicit {R}unge-{K}utta implementations with
  optimal multilevel preconditioners at the scaling limit}.
\bjtitle{SIAM J. Sci. Comput.}
\bvolume{46}(\bissue{2}),
\bfpage{71}--\blpage{96}
(\byear{2024})
\doiurl{10.1137/22M1503270}
\end{barticle}
\endbibitem

\bibitem[\protect\citeauthoryear{Saad}{2003}]{SAAD}
\begin{bbook}
\bauthor{\bsnm{Saad}, \binits{Y.}}:
\bbtitle{Iterative Methods for Sparse Linear Systems},
\bedition{2}nd edn.,
p. \bfpage{528}.
\bpublisher{Society for Industrial and Applied Mathematics},
\blocation{Philadelphia, PA}
(\byear{2003}).
\doiurl{10.1137/1.9780898718003}
\end{bbook}
\endbibitem

\bibitem[\protect\citeauthoryear{Scherer and T\"{u}rke}{1983}]{SchererTurke}
\begin{barticle}
\bauthor{\bsnm{Scherer}, \binits{R.}},
\bauthor{\bsnm{T\"{u}rke}, \binits{H.}}:
\batitle{Reflected and transposed {R}unge-{K}utta methods}.
\bjtitle{BIT}
\bvolume{23}(\bissue{2}),
\bfpage{262}--\blpage{266}
(\byear{1983})
\doiurl{10.1007/BF02218447}
\end{barticle}
\endbibitem

\bibitem[\protect\citeauthoryear{Simoncini}{2007}]{MR2318706}
\begin{barticle}
\bauthor{\bsnm{Simoncini}, \binits{V.}}:
\batitle{A new iterative method for solving large-scale {L}yapunov matrix
  equations}.
\bjtitle{SIAM J. Sci. Comput.}
\bvolume{29}(\bissue{3}),
\bfpage{1268}--\blpage{1288}
(\byear{2007})
\doiurl{10.1137/06066120X}
\end{barticle}
\endbibitem

\bibitem[\protect\citeauthoryear{Simoncini}{2016}]{SimonciniReview}
\begin{barticle}
\bauthor{\bsnm{Simoncini}, \binits{V.}}:
\batitle{Computational methods for linear matrix equations}.
\bjtitle{SIAM Rev.}
\bvolume{58}(\bissue{3}),
\bfpage{377}--\blpage{441}
(\byear{2016})
\doiurl{10.1137/130912839}
\end{barticle}
\endbibitem

\bibitem[\protect\citeauthoryear{Sommeijer}{1993}]{Sommeijer1}
\begin{barticle}
\bauthor{\bsnm{Sommeijer}, \binits{B.P.}}:
\batitle{Explicit, high-order {R}unge-{K}utta-{N}ystr\"{o}m methods for
  parallel computers}.
\bjtitle{Appl. Numer. Math.}
\bvolume{13}(\bissue{1-3}),
\bfpage{221}--\blpage{240}
(\byear{1993})
\doiurl{10.1016/0168-9274(93)90145-H}
\end{barticle}
\endbibitem

\bibitem[\protect\citeauthoryear{Southworth et~al.}{2022a}]{SouthworthII}
\begin{barticle}
\bauthor{\bsnm{Southworth}, \binits{B.S.}},
\bauthor{\bsnm{Krzysik}, \binits{O.}},
\bauthor{\bsnm{Pazner}, \binits{W.}}:
\batitle{Fast solution of fully implicit {R}unge-{K}utta and discontinuous
  {G}alerkin in time for numerical {PDE}s, {P}art {II}: {N}onlinearities and
  {DAE}s}.
\bjtitle{SIAM J. Sci. Comput.}
\bvolume{44}(\bissue{2}),
\bfpage{636}--\blpage{663}
(\byear{2022})
\doiurl{10.1137/21M1390438}
\end{barticle}
\endbibitem

\bibitem[\protect\citeauthoryear{Southworth et~al.}{2022b}]{SouthworthI}
\begin{barticle}
\bauthor{\bsnm{Southworth}, \binits{B.S.}},
\bauthor{\bsnm{Krzysik}, \binits{O.}},
\bauthor{\bsnm{Pazner}, \binits{W.}},
\bauthor{\bsnm{De~Sterck}, \binits{H.}}:
\batitle{Fast solution of fully implicit {R}unge-{K}utta and discontinuous
  {G}alerkin in time for numerical {PDE}s, {P}art {I}: {T}he linear setting}.
\bjtitle{SIAM J. Sci. Comput.}
\bvolume{44}(\bissue{1}),
\bfpage{416}--\blpage{443}
(\byear{2022})
\doiurl{10.1137/21M1389742}
\end{barticle}
\endbibitem

\bibitem[\protect\citeauthoryear{van~der Houwen and Sommeijer}{1991}]{Houwen1}
\begin{barticle}
\bauthor{\bsnm{Houwen}, \binits{P.J.}},
\bauthor{\bsnm{Sommeijer}, \binits{B.P.}}:
\batitle{Iterated {R}unge-{K}utta methods on parallel computers}.
\bjtitle{SIAM J. Sci. Statist. Comput.}
\bvolume{12}(\bissue{5}),
\bfpage{1000}--\blpage{1028}
(\byear{1991})
\doiurl{10.1137/0912054}
\end{barticle}
\endbibitem

\bibitem[\protect\citeauthoryear{van~der Houwen and Sommeijer}{1994}]{Houwen3}
\begin{botherref}
\oauthor{\bsnm{Houwen}, \binits{P.J.}},
\oauthor{\bsnm{Sommeijer}, \binits{B.P.}}:
Preconditioning in parallel {R}unge-{K}utta methods for stiff initial value
  problems.
vol. 28,
pp. 17--31
(1994).
\doiurl{10.1016/0898-1221(94)00183-9}
\end{botherref}
\endbibitem

\bibitem[\protect\citeauthoryear{van~der Houwen et~al.}{1992}]{Houwen2}
\begin{barticle}
\bauthor{\bsnm{Houwen}, \binits{P.J.}},
\bauthor{\bsnm{Sommeijer}, \binits{B.P.}},
\bauthor{\bsnm{Couzy}, \binits{W.}}:
\batitle{Embedded diagonally implicit {R}unge-{K}utta algorithms on parallel
  computers}.
\bjtitle{Math. Comp.}
\bvolume{58}(\bissue{197}),
\bfpage{135}--\blpage{159}
(\byear{1992})
\doiurl{10.2307/2153025}
\end{barticle}
\endbibitem

\bibitem[\protect\citeauthoryear{Verdugo and Badia}{2022}]{Verdugo2022}
\begin{barticle}
\bauthor{\bsnm{Verdugo}, \binits{F.}},
\bauthor{\bsnm{Badia}, \binits{S.}}:
\batitle{The software design of {G}ridap: {A} {F}inite {E}lement package based
  on the {J}ulia {JIT} compiler}.
\bjtitle{Computer Physics Communications}
\bvolume{276},
\bfpage{108341}
(\byear{2022})
\doiurl{10.1016/j.cpc.2022.108341}
\end{barticle}
\endbibitem

\end{thebibliography}
